\documentclass[a4paper]{amsart} 

\usepackage[colorlinks,linkcolor=blue,citecolor=blue,urlcolor=blue]{hyperref}
\usepackage[english]{babel}
\usepackage{amsthm,amsfonts,amssymb,mathrsfs,amsmath}
\usepackage{enumitem}
\usepackage{xcolor}
\usepackage{soul}
\usepackage{algorithm}
\usepackage{algpseudocode}

\usepackage{mathtools}						
\mathtoolsset{showonlyrefs,showmanualtags}	

\usepackage{todonotes}  

\newtheorem{Theorem}{Theorem}[section]
\newtheorem{Coro}[Theorem]{Corollary}
\newtheorem{Lemma}[Theorem]{Lemma}
\newtheorem{Prop}[Theorem]{Proposition}

\theoremstyle{definition}
\newtheorem{Def}[Theorem]{Definition}

\theoremstyle{remark}
\newtheorem{rk}[Theorem]{Remark}

\newcounter{numcond}

\newtheorem{assumption}[numcond]{Assumption}

\newcommand{\eps}{\varepsilon}

\newcommand{\hh}{{\mathcal H}}

\newcommand{\Sym}{\textnormal{Sym}}

\newcommand{\R}{\mathbb{R}}

\newcommand{\N}{\mathbb{N}}
\newcommand{\C}{\mathbb{C}}

\numberwithin{equation}{section}

\title[Inverse non autonomous linear-quadratic problem]{Inverse optimal control problem in the non~autonomous linear-quadratic case}

\author{Fr\'ed\'eric Jean}
\address{Fr\'ed\'eric Jean \\
	UMA, ENSTA Paris, Institut Polytechnique de Paris, F-91120 Palaiseau, France}
\email{frederic.jean@ensta-paris.fr}
\urladdr{\url{http://uma.ensta-paris.fr/~fjean}}

\author{Sofya Maslovskaya}
\address{Sofya Maslovskaya \\
	Department of Mathematics, Paderborn University, Paderborn, Germany}
\email{sofya.maslovskaya@uni-paderborn.de}

\begin{document}
\begin{abstract}		
Inverse optimal control problem emerges in different practical applications, where the goal is to design a cost function in order to approximate given optimal strategies of an expert. Typical application is in robotics for generation of human motions. In this paper we analyze a general class of non autonomous inverse linear quadratic problems. This class of problems is of particular interest because it arises as a linearization of a nonlinear problem around an optimal trajectory. The addressed questions are the injectivity of the inverse problem and the reconstruction. We show that the nonlinear problem admits the same characterization of the injectivity as the autonomous one. In the autonomous case we show moreover that the injectivity property is generic in the considered class. We also provide a numerical test of the reconstruction algorithm in the autonomous setting.
\end{abstract}
\maketitle


\section{Introduction}
The movements of living creatures \cite{McFarland1977,Alexander1997}, including humans \cite{Arechavaleta2008,Berret2011}, can be described by optimality principles. In particular, optimal control setting proved to be efficient in modeling of human movements \cite{Berret2008,Chitour2012} and other organisms \cite{Bettiol2019,Giraldi2020}. In most cases the optimality criteria governing the movements are unknown, which makes it an important task to find such criteria. This leads us to the inverse optimal control problem, where the goal is to recover optimality principles underlying the observed movements. Such problem is especially important for medical purposes, where it can be used for understanding of the neurophysiological behavior, designing
a rehabilitation exoskeleton well adapted to the patient, or microrobots for delivery of medical substances. It is also a substantial problem in robotics and contributes to robot
movements inspired by biological ones, autonomous vehicles and robots moving and acting in coordination with humans.

The mathematical formulation of the inverse optimal control problem consists in finding a cost function such that the given set of trajectories is minimizing for this cost under the known control system constraints. In this paper we study the inverse problem for the class of time-varying linear quadratic problems (inverse LQ problem for short). The time-varying setting is of high interest in the inverse optimal control problems arising in robotics \cite{Westermann2020}. In the autonomous framework, the problem was first stated and treated by Kalman for the single input case in infinite horizon setting \cite{Kalman1964}. Then, it was generalized to the multiinput case \cite{Anderson1966}. And for finite horizon it was treated in \cite{Nori2004,Jean2018s,Li2020}. The only instance of the treatment of the time-varying linear quadratic problem is in \cite{Jameson1973} for the best of our knowledge. The authors of \cite{Jameson1973} treat the inverse problem for a class of time-varying problems with free final condition and the given data are optimal control strategies in feedback form. We extend these results by treating the problem with fixed initial and final conditions and the input of the inverse problem given by the set of optimal trajectories in the state space. One of the important application of the time-varying linear quadratic problem is the linearization of a nonlinear optimal control problem around a trajectory. Notice that the linearized control system around periodic strokes are time-periodic. All “uniformity” hypotheses used in the assumptions of this paper are satisfied for the LQ inverse problem arising from the linearisation of the original one around periodic stroke. 

One of the fundamental questions in the inverse optimal control problem is the well-posedness, where the main issues are the existence, uniqueness and stability of the solution. The existence can be usually assumed from the knowledge of the physical system, so that the next question to treat is the uniqueness. It is possible to set the inverse problem in a more abstract setting by defining a function which maps a cost function to the corresponding solutions of the optimal control problem with fixed control system. The uniqueness of solutions to the inverse problem translates into the injectivity of the defined function. The injectivity does not hold in general even if we disregard the simple case of constantly proportional costs by considering some kind of normalization. In the series of works \cite{Jean2018a,Jean2019a,Jean2019}, this question was treated for some classes of affine control systems with quadratic cost including autonomous linear quadratic problem. In this work we extend these results to the non autonomous LQ problem. 

In the direct optimal control problems, the necessary optimality conditions for a trajectory $x(\cdot)$ is that it can be lifted to a trajectory $(x(\cdot), p(\cdot))$, called extremal, of a Hamiltonian system defined in the cotangent bundle of the state space. Projections of extremals to the state space are called geodesics and they are candidates to be solution to the direct optimal control problem. Inverse optimal control problem can be characterised in a rather general non-linear setting (control-affine systems with quadratic costs \cite{Jean2019a} or more generally with Tonelli Lagrangian \cite{Maslovskaya2018}) as follows. If two costs have the same minimizers, then, locally near ample geodesics (projections to the state space of a particular type of extremals), there exists an orbital diffeomorphism between the corresponding Hamiltonian vector fields (\cite[Prop. 3.6]{Jean2019a} and \cite[Prop. 1.41]{Maslovskaya2018}).
Note that
\begin{itemize}
	\item the linearization of the Hamiltonian extremal flow around an extremal $(x(\cdot), p(\cdot))$ is a linear Hamiltonian flow associated with the linearized control system around $x(\cdot)$ and with a quadratic cost;
	
	\item by definition, a geodesic of a control-affine system is ample if the linearized 
	system satisfies Kalman sufficient condition (I.3).
\end{itemize}
By differentiation of the relation given by the orbital diffeomorphism, one gets a one-to-one correspondence between the trajectory of the linearized Hamiltonian and the two LQ problems associated with the linearized control system having the
same minimizers. More details on the connection between the the nonlinear inverse problem and the corresponding linearization can be found in the Appendix.

In this paper we generalise all the results obtained in the autonomous LQ problem \cite{Jean2018a} to the non autonomous case and we show that a similar description of the costs leading to the same solution to optimal control problem can be obtained. This is done in several steps.
\begin{enumerate}
	\item We find a form of the optimal solution similar to the one in the autonomous case \cite{Jean2018a} which permits to parametrise the optimal synthesis by pairs of matrix-valued functions.
	\item We formulate a class of canonical costs similar to the canonical costs in \cite{Jean2018a} for which we provide the injectivity criteria. The non-injective cases are characterised by the product structure, which is in line with the results obtained in \cite{Jean2018a}, but this time in the non-autonomous context. Moreover, we show the genericity of the injective case for the autonomous LQ problem.
	\item The reconstruction is discussed for the general non-autonomous inverse LQ problem and a concrete algorithm is given for the autonomous case, for which we also present the numerical results showing its stability.
\end{enumerate}

The paper is organized as follows. In Section~\ref{sec:non.autonome} we characterize the solutions of the direct LQ problem in the non autonomous setting and formalize the inverse problem. The main results on the injectivity are contained in Section~\ref{se:injectivity} as well as a discussion on the reconstruction approach. Section~\ref{sec:autonomous} presents the new results in the autonomous case including the genericity of the injective case and numerical results of the reconstruction algorithm.
%
%
%
%
%
%
%
%
%
%
\section{Direct and inverse linear quadratic optimal control problems}
\label{sec:non.autonome}
Let us consider a time-dependent linear control system
\begin{equation}
\label{eq:dymsys}
\dot{x}(t) = A(t)x(t) + B(t)u(t), \qquad x(t)\in \R^n, \quad u(t) \in \R^m, \quad t \in \R,
\end{equation}
and a quadratic cost of the form
\begin{equation}
\label{eq:qvcost}
c(t,x,u)=  x^{\top}Q(t) x + 2 x^{\top} S(t) u + u^{\top} R(t) u,
\end{equation}
where $A(t),B(t),Q(t),S(t),R(t)$ are real-valued matrices of the appropriate dimensions. These control system and cost
define a family of linear-quadratic (LQ) optimal control problems with terminal constraints: given a pair of times $t_0<t_1$, an initial point $x_0 \in \R^n$ and a final point $x_1\in \R^n$, minimize the integral cost
\begin{equation}
C_{t_0}^{t_1}(x_u) = \int_{t_0}^{t_1} c(t,x_u(t),u(t)) dt
\end{equation}
among all trajectories $x_u$ of \eqref{eq:dymsys} associated with the control law $u$ satisfying $x_u(t_0) = x_0$ and $x_u(t_1) = x_1$.

Under suitable assumptions (see Theorem~\ref{th:decomp_sol} below) on the mappings $A,B,Q$, $S,R$, every LQ problem in the above family admits a unique minimizing trajectory on any interval $[t_0,t_1]$. The set of all the minimizing solutions, for all data $t_0,t_1,x_0,x_1$, is called the \emph{optimal synthesis} associated with the dynamics~\eqref{eq:dymsys} and the cost~\eqref{eq:qvcost} (we say also that the synthesis is associated with $(A,B,Q,S,R)$).

An \emph{inverse linear-quadratic optimal control problem} is posed as follows: a linear system \eqref{eq:dymsys} being fixed, recover the quadratic cost $c$ from the given of the corresponding optimal synthesis. \medskip

To study such an inverse problem,
it is necessary to discuss its well-posedness and in particular its injectivity: is there a unique cost $c$ associated with a given optimal synthesis? The aim of this section is to fix a framework in which this question  may have a satisfactory answer, we then characterize the uniqueness of solutions in the next one, together with the reconstruction of the cost. Let us first introduce the characterization of optimal syntheses that we will use.


\subsection{Characterization of the optimal synthesis}
Let us recall the two usual ways to describe the solution of a LQ problem with terminal constraints. The first one arises from the Pontryagin Maximum Principle. Given data $t_0,t_1,x_0,x_1$, the trajectory solution of the corresponding LQ problem is defined by the control law
\begin{equation}
u(t)= - R(t)^{-1}S(t)^\top x(t) + R(t)^{-1}B(t)^\top p(t), \quad t \in [t_0,t_1],
\end{equation}
where $(x,p)(\cdot)$ is the unique solution of the Hamiltonian differential equation
\begin{equation}
\label{eq:ham_edo}
\begin{pmatrix}
    \dot x (t) \\
    \dot p (t) \\
\end{pmatrix} = \vec{H}(t)
\begin{pmatrix}
     x (t) \\
     p (t) \\
\end{pmatrix},  \quad t \in [t_0,t_1],
\end{equation}
such that $x(t_0) = x_0$ and $x(t_1) = x_1$, the Hamiltonian matrix being defined by
\begin{equation}
\label{eq:Ham_matrix}
\vec{H}(t) =  \begin{pmatrix}
		A(t) - B(t)R(t)^{-1}S(t)^\top & B(t)R(t)^{-1}B(t)^\top\\
		Q(t) - S(t)R(t)^{-1}S(t)^\top & -A(t)^\top + S(t)R(t)^{-1}B(t)^\top
		\end{pmatrix}.
\end{equation}
The drawback of this description is that the computation of $p(t)$  depends explicitly on all data, so we need the collection of all possible functions $p$ to describe the optimal synthesis.

The second way of describing solutions arises through the resolution of the Riccati differential equation (we omit the time dependance for readability),
\begin{equation}\label{eq:Riccati_edo}
\dot P =  - PA - A^\top P + (PB + S) R^{-1} (B^\top P + S^\top) -Q.
\end{equation}
Indeed (see for instance~\cite{Bryson1975}), the optimal control law is given by
\begin{equation}
u(t)= - R(t)^{-1}(B(t)^\top M(t)+S(t)^\top) x(t) - R(t)^{-1}B(t)^\top N(t) x_1, \quad t \in [t_0,t_1],
\end{equation}
where $M(t)=P(t)-V(t)W(t)^{-1}V(t)^\top$ and $N(t)=V(t)W(t)^{-1}$, the matrix-valued functions $P,V,W$ being defined as follows: $P$ is the solution on $[t_0,t_1]$ of the Riccati differential equation~\eqref{eq:Riccati_edo} such that $P(t_1)=0$, and $V,W$ are the solutions of
\begin{align}
  \dot V & =  \left( -A^\top + (PB+S)R^{-1} B^\top \right) V, \qquad V(t_1)=I, \\
  \dot W & =  V^\top B R^{-1}B^\top V, \qquad W(t_1)=0.
\end{align}
This description of $u$ as a feedback $u(t)=K(t)x(t)+K_1(t)x_1$ is more practical since the computation of the gains $K,K_1$ allows to determine all the optimal controls for all values of $x_0,x_1$.  However these gains still depend explicitly on the terminal time $t_1$, so we need the collection of all the gains  $K,K_1$ for all $t_1 \in \R$ in order to determine the optimal synthesis. Again, this characterization is unusable in the context of the inverse optimal control.

We propose here a third way of describing optimal solutions that allows us to characterize easily the optimal synthesis. This description  requires however stronger assumptions than what is usually required. We detail these assumptions now, distinguishing between those on dynamics and those on cost. Let us first introduce some notations. We use $M_n(\R)$ and $M_{n,m}(\R)$ to denote the set of real-valued $(n \times n)$ and $(n \times m)$ matrices, respectively, $GL_n(\R)$ to denote the set of invertible matrices in $M_n(\R)$, and $\Sym(n)$ to denote the set of symmetric matrices in $M_n(\R)$. For $Q_1, Q_2 \in \Sym(n)$, we write $Q_1 \succ Q_2$ if $Q_1-Q_2$ is positive definite. Moreover,  we use $\Phi_A(t,\tau)$ to denote the \emph{transition matrix} of the differential equation $\dot x (t)=A(t) x(t)$, i.e., $\Phi_A(\cdot,\tau)$ is the solution of
\begin{equation}
   \frac{\partial}{\partial t}\Phi_A(t,\tau) = A(t) \Phi_A(t,\tau), \qquad \Phi_A(\tau,\tau)=I.
\end{equation}
Finally, the mapping $A$ is said to be \emph{exponentially stable} is there exists two positive constants $\alpha,\beta$ such that
\begin{equation}
  \|\Phi_A (t,t_0) \| \leq \alpha e^{-\beta (t-t_0)} \quad \hbox{for every } t \geq t_0,
\end{equation}
and \emph{exponentially anti-stable} is there exists two positive constants $\alpha,\beta$ such that
\begin{equation}
  \|\Phi_A (t,t_0) \| \leq \alpha e^{\beta (t-t_0)} \quad \hbox{for every } t \leq t_0.
\end{equation}

The assumption on the dynamics $\dot x = A x + Bu$ are the following.
\begin{assumption}
	\label{Assump:AB}
	The mappings $A,B$ defining the dynamics \eqref{eq:dymsys} satisfy:
	\begin{enumerate}[label=\textit{(\ref{Assump:AB}.\arabic*)}]
    \item  \label{cond:A1}
    $A : \R \to M_n(\R)$ and $B : \R \to M_{n,m}(\R)$ are smooth ($C^\infty$) and bounded mappings;
		
    \item  \label{cond:B_injectif}
    for every $t \in \R$, $B(t)$ is injective (in particular $m \leq n$); 

	\item  \label{cond:kalman}
the pair $(A,B)$ satisfies the \emph{Kalman sufficient controllability condition} everywhere, i.e.,
		\begin{equation}
        \mathrm{span}\{B_i(t)v \ : \ v \in \R^m, \ i \in \N\} = \R^n, \qquad \hbox{for every } t \in \R,
		\end{equation}
where the smooth mappings $B_i$ are defined inductively by $B_0(t)=B(t)$, $B_{i+1}(t)=\dot{B}_i(t) - A(t)B_i(t)$;

\item \label{cond:unif_contr}
the pair $(A ,B)$ is \emph{uniformly controllable}: there exist positive constants $s, \beta >0$  such that, for every $t \in \R$,
\begin{equation}\label{eq:controllability.matrix}
  \int_{t}^{t+s} \Phi_{A}(t,\tau) B(\tau) B(\tau)^\top \Phi_{A}(t,\tau)^\top d\tau \succ \beta I.
  \end{equation}

\end{enumerate}
\end{assumption}

\begin{rk}
\label{re:assumptions}
Different notions of controllability appear in these assumptions, let us recall their meaning (we refer to classical text books such as~\cite{Anderson1990} and~\cite[Chap.\ 1]{Coron2007} for the results below).
\begin{enumerate}
  \item Uniform controllability (condition \ref{cond:unif_contr}) implies not only that the system can be steered from
any state to any other one, but also that the minimal $L^2$ norm of the control needed in this
transfer and the transfer time are roughly independent of the initial time. This property is preserved by bounded feedback: if \ref{cond:unif_contr} holds, then for any bounded mapping $K : \R \to M_{m,n}(\R)$, the pair $(A-BK,B)$ is also uniformly continuous.  Moreover, \ref{cond:unif_contr} implies that there exists bounded mappings $K_+,K_-$ such that $A-BK_+$ is exponentially stable and $A-BK_-$ is exponentially anti-stable.

  \item \ref{cond:kalman} and \ref{cond:unif_contr}  are independent while both imply controllability. It is for instance easy to construct examples of pairs $(A,B)$ which satisfy \ref{cond:kalman} but are not uniformly controllable: take for instance $A(t)=A_0$ and $B(t)=e^{-t}B_0$, where $(A_0,B_0)$ are constant matrices satisfying the Kalman rank condition. Conversely there exist non analytic smooth mappings $A,B$ which are controllable (and even uniformly controllable) but which nonetheless do not satisfy \ref{cond:kalman}.

  \item A straightforward induction argument shows that \ref{cond:kalman} implies that any pair of the form $(A-BK,B)$ or $(A,BU)$, where $K: \R \to M_{m,n}(\R)$ and $U: \R \to GL_m(\R)$ are smooth, also satisfies Kalman sufficient controllability condition. 
\end{enumerate}
\end{rk}

\begin{rk}
	\label{re:periodic}
In the case where the mappings $A$ and $B$ are periodic with the same period $T>0$, Kalman condition \ref{cond:kalman} implies uniform controllability~\ref{cond:unif_contr}. Indeed, Kalman condition \ref{cond:kalman} implies that the system is controllable on any time interval $[t_0, t_1], \ t_0 < t_1 \in \R$. Hence, by \cite[Th.\ 1.11]{Coron2007},  the Gramian satisfies for every $t \in \R$
$$
\int_{t}^{t+T}  \Phi_{A}(t,\tau) B(\tau) B(\tau)^\top \Phi_{A}(t,\tau)^\top \succ \beta(t) I,
$$
where $\beta$ is a positive continuous function. 
	If $A$ and $B$ are  $T$-periodic, then the transition matrix of $A$ satisfies $\Phi_A(t,\tau)=\Phi_A(t+T,\tau + T)$ and we have
\begin{multline}
\int_{t}^{t+T} \Phi_{A}(t,\tau) B(\tau) B(\tau)^\top \Phi_{A}(t,\tau)^\top d\tau = \\
\int_{t+T}^{t+2T} \Phi_{A}(t+T,\tau) B(\tau) B(\tau)^\top \Phi_{A}(t+T,\tau)^\top d\tau.
\end{multline}
As a consequence $\beta$ is $T$-periodic and then admits a positive minimum, which implies the uniform controllability~\ref{cond:unif_contr}. Note that both conditions \ref{cond:kalman} and \ref{cond:unif_contr} coincide if $A,B$ are moreover analytic functions, because \ref{cond:kalman} is also a necessary condition of controllability in this case.
\end{rk}

\begin{rk}
	\label{re:periodic.linearization}
	An example of a control system which satisfy Assumptions~\ref{Assump:AB} is the linearization of an affine control system around a periodic ample geodesic. In this case, \ref{cond:kalman} is satisfied by definition of ample geodesics \cite[Prop.\ 3.12]{AgrachevBarilariRizzi}, and on the other hand, linearization around a periodic trajectory is periodic, which implies also \ref{cond:unif_contr} by the preceding remark. As it was noted in introduction, periodic trajectories of affine control systems are of particular interest in applications of inverse optimal control.
\end{rk}

The assumptions on the quadratic cost are the following ones.
\begin{assumption}
	\label{Assump:ABQSR}
The mappings $(Q,S,R)$ defining the cost $c$ satisfy:
\begin{enumerate}[label=\textit{(\ref{Assump:ABQSR}.\arabic*)}]
\item \label{cond:B1}
$Q : \R \to \Sym(n)$,
 $S : \R \to M_{n,m}(\R)$ and $R : \R \to \Sym(m)$ are smooth and bounded mappings;

  \item \label{cond:Rdefpos}
  $R$ is \emph{uniformly positive definite}: $\exists \alpha >0$ such that $R(t) \succ \alpha I$ for any $t \in \R$;

\item \label{cond:barQpos}
for every $t \in \R$, the symmetric matrix $Q_S(t)= Q(t) - S(t) R(t)^{-1} S(t)^\top$ is nonnegative.


%
\end{enumerate}
\end{assumption}

\begin{rk}
\label{re:assumptions2}
The conditions appearing in Assumptions~\ref{Assump:AB} and~\ref{Assump:ABQSR} play different roles.
\begin{enumerate}
  \item The boundedness hypothesis in \ref{cond:A1} and \ref{cond:B1} together with the uniform properties \ref{cond:unif_contr} and \ref{cond:Rdefpos}, and hypothesis~\ref{cond:barQpos} are necessary to apply Theorems 3.5 and 3.6 of \cite{Kwakernaak1972}, which are the basic tools to establish Lemma~\ref{le:P+} and Theorem~\ref{th:decomp_sol}.
  \item The proof of Theorem~\ref{th:decomp_sol} do not require Assumptions \ref{cond:B_injectif} and \ref{cond:kalman}, nor the smoothness hypothesis in \ref{cond:A1} and \ref{cond:B1} (continuity would be sufficient). However the first two  will be used in the study of injectivity of the inverse problem, in particular in the proof of Lemma~\ref{prop:orb.diffeo}, which is a crucial step. And the smoothness is required to write the Kalman condition \ref{cond:kalman} for $(A,B)$ and for the pairs $(A-BK,B)$ that appear at different places in the proofs.

    \item When the mappings $A,B,Q,S,R$ are only smooth (and not necessarily bounded), existence and uniqueness of the optimal controls is guaranteed by hypotheses~\ref{cond:kalman}, \ref{cond:Rdefpos}, and \ref{cond:barQpos} (see for instance~\cite[Sect.\ 3.3]{Lee1967}).

\end{enumerate}
\end{rk}

We are now in a position to state our result on the description of the optimal synthesis.

\begin{Theorem}
\label{th:decomp_sol}
Consider mappings $(A,B,Q,S,R)$ satisfying Assumptions~\ref{Assump:AB} and~\ref{Assump:ABQSR}.

\begin{enumerate}[label={(\roman*)}]
  \item There exists a unique solution $P_+ : \R \to \Sym(n)$ of the Riccati differential equation~\eqref{eq:Riccati_edo} such that the matrix-valued mapping
  \begin{equation}
\label{eq:A+}
A_{+} = A-B R^{-1} \left(S^\top + B^\top P_{+}\right)
\end{equation}
is exponentially stable. Moreover, for every $t \in \R$, $P_+(t)$ is nonnegative.

\item There exists a unique solution $P_- : \R \to \Sym(n)$ of the Riccati differential equation~\eqref{eq:Riccati_edo} such that the matrix-valued mapping
\begin{equation}
\label{eq:A-}
A_{-} = A-B R^{-1} \left(S^\top + B^\top P_{-}\right)
\end{equation}
is exponentially anti-stable.
Moreover, for every $t \in \R$, $P_-(t)$ is nonpositive.

\item For every $t \in \R$, $P_+(t) - P_-(t)$ is positive definite.

\item \label{it:inutile} For every data $t_0<t_1$, $x_0,x_1 \in \R^n$, there exists a unique trajectory $x(\cdot)$ of the corresponding LQ optimal problem associated with $(A,B,Q,S,R)$. This solution is the first component of the unique solution $(x,p)$ of the Hamiltonian differential equation \eqref{eq:ham_edo} such that $x(t_0) = x_0$ and $x(t_1) = x_1$. Moreover there exists a unique pair $(y_+,y_-) \in \R^{2n}$ such that
\begin{equation}
\label{eq:decomp_sol}
x(t) = \Phi_{A_+}(t,t_0) y_+ + \Phi_{A_-}(t,t_0) y_-, \qquad \hbox{for every } t \in [t_0,t_1].
\end{equation}
\end{enumerate}
\end{Theorem}

This result implies that an optimal synthesis is completely characterized by the pair of mappings $(A_+,A_-)$. Note that the description \eqref{eq:decomp_sol} is already known in the autonomous case (see~\cite{Ferrante2005}), but to the best of our knowledge it is new in the non autonomous case. The proof of the theorem is quite long and does not fall within the main scope of the article, so it is postponed  to  Appendix.~\ref{se:appendix}.


\subsection{Injectivity and cost equivalence}
\label{se:reduc}
Fix a pair of mappings $(A, B)$  satisfying Assumptions \ref{Assump:AB} and denote by $\mathcal{C}_{LQ}$ the set of quadratic costs $c$ of the form \eqref{eq:qvcost} satisfying Assumptions \ref{Assump:ABQSR}.

We can not expect that the map which to a cost associates the corresponding optimal synthesis is injective on  the set $\mathcal{C}_{LQ}$. Indeed, it is easy to construct different costs having the same optimal synthesis, for instance by taking two costs constantly proportional to each other (we will see other examples in next subsection). This leads us to the following definition.
\begin{Def} \label{def:equiv}
We say that two costs $c, \tilde{c}$ in  $\mathcal{C}_{LQ}$ are \emph{equivalent}, and we write $c \sim \tilde{c}$, if they define the same optimal synthesis.
\end{Def}

Theorem~\ref{th:decomp_sol} directly implies  that an optimal synthesis of a cost in $\mathcal{C}_{LQ}$ is completely characterized by the pair of mappings $(A_+,A_-)$. We prove now that this characterization is univocal.

\begin{Lemma}
\label{lem:A+A-}
Two equivalent costs define the same mappings $(A_+, A_-)$.
\end{Lemma}
\begin{proof}
	Consider two equivalent costs  $c$ and $\tilde{c}$ defined by the mappings $(Q, S, R)$ and $(\tilde{Q}, \tilde{S}, \tilde{R})$, respectively. By definition of equivalence, the two costs define the same minimizing trajectories.
	
	Fix $t_0 \in \R$. For $t_1>t_0$ and $i=1,\dots,n$, let $x_i(\cdot)$ be the minimizing trajectory between $x_i(t_0)=e_i$ and $x_i(t_1)=\Phi_{\tilde{A}_+}(t_1,t_0)e_i$, where $e_i$ denotes the $i$th vector of the canonical basis of $\R^n$. By uniqueness of the decomposition~\eqref{eq:decomp_sol} for $\tilde{c}$, we get $x_i(t) = \Phi_{\tilde{A}_+}(t,t_0)e_i$. We write this in matrix form  as
\begin{equation}
X(t)=
\begin{pmatrix}
x_1(t) & \cdots & x_n(t)
\end{pmatrix} = \Phi_{\tilde{A}_+}(t,t_0), \qquad t \in [t_0,t_1].
\end{equation}
	
Now, since $c \sim \tilde{c}$, there exist $(n \times n)$ matrices $Y_+^{t_1},Y_-^{t_1}$ such that
$X(t) = \Phi_{\tilde{A}_+}(t,t_0) = \Phi_{A_+}(t,t_0) Y_+^{t_1} + \Phi_{A_-}(t,t_0) Y_-^{t_1}$ for every $t \in [t_0,t_1]$. Notice that $Y_{\pm}^{t_1}$ are actually independent of $t_1$. Indeed, given $t_1 \neq t'_1$ greater than $t_0$, we have, for every $t \in [t_0, \min(t_1,t'_1)]$,
\begin{multline}
\Phi_{A_+}(t,t_0) Y_+^{t_1} + \Phi_{A_-}(t,t_0) Y_-^{t_1} =
\Phi_{\tilde{A}_+}(t,t_0) = \Phi_{A_+}(t,t_0) Y_+^{t'_1} + \Phi_{A_-}(t,t_0) Y_-^{t'_1},
\end{multline}
which by uniqueness of the decomposition~\eqref{eq:decomp_sol} implies $Y_{\pm}^{t_1}=Y_{\pm}^{t'_1}$. As a consequence, setting $Y_{\pm}:=Y_{\pm}^{t_1}$, there holds
\begin{equation}
\label{eq:PhiA-tilde}
\Phi_{A_+}(t,t_0) Y_+ + \Phi_{A_-}(t,t_0) Y_- =
\Phi_{\tilde{A}_+}(t,t_0) \quad \text{for any} \ t \in [t_0, +\infty).
\end{equation}
As $t \to +\infty$, $\Phi_{A_+}(t,t_0)$ and $\Phi_{\tilde{A}_+}(t,t_0)$ tend to zero since $A_+$ and $\tilde{A}_+$ are asymptotically stable, hence $\Phi_{A_-}(t,t_0)Y_-  \to 0$. This in turn implies $Y_-=0$ since $A_-$ is exponentially anti-stable, and then $\Phi_{A_+}(t,t_0) Y_+  = \Phi_{\tilde{A}_+}(t,t_0)$ for any $t \geq t_0$. Taking $t=t_0$ we obtain  $Y_+=I$, and differentiating the equality with respect to $t$ we get $A_+ (t)= \tilde{A}_+(t)$ for $t \geq t_0$.

The time $t_0$ having been fixed arbitrarily, we finally get $A_+ = \tilde{A}_+$. By exchanging the role of $A_+$ and $A_-$ and taking $t \to - \infty$, we obtain in the same way that $ A_- = \tilde{A}_-$.
\end{proof}

To summarize, we saw in Theorem~\ref{th:decomp_sol} that a LQ optimal synthesis is completely characterized by the pair of mappings $(A_+, A_-)$, and by the previous lemma we showed that this pair is unique. In other terms, to every equivalence class in $\mathcal{C}_{LQ}$ corresponds a unique pair $(A_+, A_-)$. Note also that, the pair $(A,B)$ being fixed, Assumption~\ref{cond:B_injectif} implies that the knowledge of $(A_+, A_-)$ is equivalent to the one of the mappings $K_+,K_-: \R \to M_{n,m}(\R)$ defined by
\begin{equation}
\label{eq:A+A-viaK+K-}
K_{+} = R^{-1} \left(S^\top + B^\top P_{+}\right) \quad \text{and} \quad K_{-} = R^{-1} \left(S^\top + B^\top P_{-}\right),
\end{equation}
since $A_{+} = A-BK_{+}$ and $A_{-} = A-BK_{-}$.

\subsection{Inverse problem in canonical form}	
\label{se:reduction}
To address the problem of injectivity, we will reduce the inverse problem to a special category of canonical costs containing a representative of each equivalence class of $\mathcal{C}_{LQ}$. The idea of restriction to some smaller category was proposed first in \cite{Nori2004} for single input autonomous problems, we will further develop the idea and construct a new category of LQ problems.

\begin{Lemma}
	\label{le:canonical}
	A quadratic cost $c \in \mathcal{C}_{LQ}$ of the form \eqref{eq:qvcost} is equivalent to the cost
	\begin{equation}
	\label{eq:cost3}
	\tilde{c} = (u+K_+x)^\top R (u+K_+x), \quad \text{with} \quad K_{+} = R^{-1} \left(S^\top + B^\top P_{+}\right),
	\end{equation}
the mapping $P_+$ being the one given by Theorem~\ref{th:decomp_sol}.
\end{Lemma}
\begin{proof}
	Given data $t_0<t_1$, $x_0, x_1$, let $x^*(\cdot)$ be the solution of $\min \int_{t_0}^{t_1} c(t,x,u) dt$ between $x_0$ and $x_1$. Clearly, $x^*(\cdot)$ minimizes as well the functional
	$$
	x_1^\top P_+ (t_1) x_1 - x_0^\top P_+ (t_0) x_0 + \int_{t_0}^{t_1} c (t, x(t),u(t)) dt.
	$$
It is a standard result on Riccati differential equation (that can be checked easily by a direct computation) that the above functional  can be written in integral form as
	$$
	\int_{t_0}^{t_1} \left( \left(u+R^{-1} \left(S^\top + B^\top P_{+}\right) x \right)^\top R \left(u+ R^{-1} \left(S^\top + B^\top P_{+}\right)x\right) \right) dt.
	$$
Therefore $x^*(\cdot)$ minimizes also $\int_{t_0}^{t_1} \tilde{c}(t,x,u) dt$.
	We conclude that any trajectory in the optimal synthesis of $c$ also belongs to the optimal synthesis of  $\tilde{c}$, which ends the proof.
\end{proof}

This result leads us to introduce the following category of quadratic costs.

\begin{Def}
	\label{def:KR}
A \emph{canonical cost} is a quadratic cost  of the form
	$$
	c(t,x,u)= (u+K(t) x)^\top R(t) (u+K(t)x),
	$$
where:
\begin{itemize}
  \item $R : \R \to \Sym(m)$ is a smooth, bounded, and uniformly positive definite mapping such that $\det R(0)=1$;
  \item $K : \R \to M_{m,n}(\R)$ is a smooth, bounded, and stabilizing matrix, i.e., $A - BK$ is exponentially stable.
\end{itemize}
\end{Def}

Any canonical cost satisfy Assumption~\ref{Assump:ABQSR}, and so belongs to the set $\mathcal{C}_{LQ}$. Indeed, written in the usual form~\eqref{eq:qvcost}, a canonical cost corresponds to the mappings $(Q,S,R)$, where $Q=K^\top R K$ and $S=K^\top R$, which satisfy \ref{cond:B1}, \ref{cond:Rdefpos}. And the mapping $Q_S=0$ satisfies \ref{cond:barQpos}. Note also that, because of the normalization $\det R(0)=1$, two non equal canonical costs can not be trivially equivalent (i.e.\ constantly proportional).

\begin{Prop}
	\label{prop:reduced}
	Any cost $c$ in $\mathcal{C}_{LQ}$ is equivalent to a canonical cost $\tilde{c}$. Moreover the mapping $K_+$ associated with $\tilde{c}$ is $K_+=K$ (equivalently, $A_+=A-BK$).
\end{Prop}

\begin{proof}
	From Lemma~\ref{le:canonical}, any cost $c$ in $\mathcal{C}_{LQ}$ is equivalent to a cost $\tilde{c} = (u+K x)^\top R (u+Kx)$, where $R$ and $K$ have the same property as a canonical cost except for $\det R(0)=1$. However, since two constantly proportional cost are equivalent and $\det R(t) >0$ for every $t$, we can assume moreover that $\det R (0)=1$, which proves the first statement of the lemma.
	
	We are left to prove that the mapping $A_+$ associated with $\tilde{c}$ is equal to $A-BK$. Fix $t_0<t_1$. For $i=1,\dots,n$, the minimizing solution between $x_i(t_0)=e_i$ and $x_i(t_1)=\Phi_{A-BK}(t_1,t_0)e_i$ is equal to $x_i(t)=\Phi_{A-BK}(t,t_0)e_i$, since the corresponding control $u_i=-Kx_i$ satisfies $\tilde{c}(x_i(t),u_i(t))\equiv 0$ and the minimizing solution is unique. Let us write in matrix form $X(t)=(x_1(t) \cdots x_n(t))= \Phi_{A-BK}(t,t_0)$, $t \in [0,T]$.  Now, from \eqref{eq:decomp_sol} there exists $(n \times n)$ matrices $Y_+,Y_-$ such that  $X(t) = \Phi_{A_+}(t,t_0) Y_+^{t_1} + \Phi_{A_-}(t,t_0) Y_-^{t_1}$ for every $t \in [t_0,t_1]$. Arguing as in the proof of Lemma~\ref{lem:A+A-} we conclude that $A_+=A-BK$, which ends the proof.
\end{proof}

\begin{rk}
Note that the set of canonical costs depends on the mappings $A,B$, while $\mathcal{C}_{LQ}$ does not. The notion of equivalence of costs depends also on $A,B$.
\end{rk}

The pair $(A,B)$ being fixed, we formulate the inverse optimal control problem among canonical costs as follows:  \emph{given a linear-quadratic optimal synthesis $\Gamma$, find a canonical cost $c$ such that $\Gamma$ is the optimal synthesis of $c$.}


Proposition~\ref{prop:reduced} ensures that this problem always has a solution, hence we concentrate now on this problem where we can expect good injectivity properties.

%

%

\section{Injective condition and reconstruction}
\label{se:injectivity}	
%

\subsection{Product structure}	

It appears that we can obtain examples of equivalent canonical costs when the problem exhibits a kind of product structures (this idea has already been used in~\cite{Jean2019}). Let us construct such an example.

Choose an integer $N>1$, and $N$ pairs of positive integers $m_i \leq n_i$, $i=1,\dots, N$. Set $m=\sum_{i=1}^N m_i$ and $n= \sum_{i=1}^N n_i$. For each $i=1,\dots, N$, choose a linear system
\begin{equation}
\dot{x}_i (t)= A_i(t) x_i(t) + B_i(t) u_i(t), \qquad t \in \R, \ x_i(t) \in \R^{n_i}, \  u_i(t) \in \R^{m_i},
\end{equation}
where $A_i,B_i$ are smooth mappings satisfying Hypothesis~\ref{cond:kalman}, and a cost $c_i(t,x_i,u_i)=u_i^\top R_i (t) u_i$ where  $R_i : \R \to \Sym(m_i)$ is a smooth mapping satisfying Hypothesis~\ref{cond:Rdefpos}. We define a linear-quadratic problem on $\R^n$ with control in $\R^m$ by setting
\begin{align}
\label{eq:prod}
& \bar{A} = \begin{pmatrix}
A_1 & & \\
& \ddots & \\
& & A_N
\end{pmatrix},  \quad \bar{B}= \begin{pmatrix}
B_1 & & \\
& \ddots & \\
& & B_N
\end{pmatrix}, \\
\hbox{and} \quad
\bar{c}=u^\top Ru & \quad \hbox{where} \quad \bar{R} = \begin{pmatrix}
R_1 & & \\
& \ddots & \\
& & R_N
\end{pmatrix}, \quad \hbox{i.e.,} \quad
\bar{c}=\sum_{i=1}^N u_i^\top R_i u_i.
\end{align}
Note that this linear-quadratic problem satisfies hypotheses~\ref{cond:kalman}, \ref{cond:Rdefpos}, and \ref{cond:barQpos}, so by Remark~\ref{re:assumptions2}--(3), existence and uniqueness of the optimal solutions are guaranteed.

Moreover, a trajectory $x(\cdot)$ is an optimal solution of the above linear-quadratic problem  if and only if $x(\cdot) = (x_1(\cdot), \dots, x_N(\cdot))$, where each $x_i(\cdot)$ is an optimal solution of the problem associated with  $A_i, B_i, R_i$. As a consequence, the cost $\bar{c}$ is equivalent to any cost of the form
\begin{equation}
\bar{c}_\lambda(t,x,u)=\sum_{i=1}^N \lambda_i u_i^\top R_i(t) u_i  , \quad \hbox{i.e.,} \quad R_\lambda(t) = \begin{pmatrix}
\lambda_1 R_1(t) & & \\
& \ddots & \\
& & \lambda_N R_N (t)
\end{pmatrix} ,
\end{equation}
where $\lambda_1,\dots, \lambda_N$ are positive real numbers.

We can extend this construction through changes of variables.
\begin{Def} \label{def:product}
We say that a LQ optimal control problem defined by a pair of mappings $(A,B)$ and a canonical cost $c=(u+Kx)^\top R (u+Kx)$ \emph{admits a product structure} if there exists an integer $N> 1$ and smooth mappings $M: \R \to GL_n(\R)$ and $U: \R \to GL_m(\R)$ such that
\begin{equation}
\bar{A} = M^{-1}(-\dot{M}+(A-BK) M), \quad \bar{B} = M^{-1}BU,  \quad \hbox{and} \quad \bar{R} =  U^\top R U
\end{equation}
are of the form \eqref{eq:prod}.
\end{Def}

\begin{Lemma}
If $(A,B)$ and $c$ admits a product structure, then $c$ admits equivalent canonical costs.
\end{Lemma}

\begin{proof}
Through the time-dependent linear change of coordinates
\begin{equation}
\bar{x}=M(t)^{-1}x, \qquad \bar{u}=U(t)^{-1}(u+K(t) x), \quad t\in \R,
\end{equation}
any trajectory of $\dot x = A x + Bu$ is transformed into a trajectory of  $\dot{\bar{x}} = \bar{A} \bar{x} + \bar{B} \bar{u}$
and the canonical cost $c(x,u)=(u+Kx)^\top R (u+Kx)$  is equal to $\bar{c}(\bar{x},\bar{u})=\bar{u}^\top \bar{R} \bar{u}$. As a consequence, any optimal trajectory of $A,B$ and $c$ is transformed into an optimal trajectory of $\bar{A}, \bar{B}$, and $\bar{c}$, and then into an optimal trajectory of $\bar{A}, \bar{B}$, and $\bar{c}_\lambda$ for any positive real numbers $\lambda_1,\dots, \lambda_N$. It then follows that $c$ is equivalent to the cost
\begin{equation}
 c_\lambda (x,u)= (u+Kx)^\top U^{-\top} \bar{R}_\lambda U^{-1} (u+Kx).
\end{equation}
Note that $\det \bar{R}_\lambda = \prod_i (\lambda_i)^{m_i} \det \bar{R}$ and $\det \bar{R}=(\det U)^2 \det R$. Therefore, if $\lambda_1,\dots, \lambda_N$ are chosen non all equal and such that $\prod_i (\lambda_i)^{m_i}=1$, then $c_\lambda$ is a canonical cost different from $c$, which proves the lemma.
\end{proof}

In view of the previous discussion, if a problem admits a product structure, then the corresponding inverse problem has many solutions. We will see in next section that the converse is true, the product structure is actually a necessary and sufficient condition for non injectivity of the inverse problem. Let us first prove that product structures are characterized by the form of their optimal syntheses.

\begin{Def}
We say that a LQ optimal synthesis characterized by mappings $A_+,A_-$ admits a \emph{product form} if there exist a non trivial decomposition $\R^n=E_1 \oplus E_2$ and a smooth mapping $M: \R \to GL_n(\R)$ such that, for every $t \in \R$, both matrices $\bar{A}_+(t)$ and $\bar{A}_-(t)$ leave invariant $E_1$ and $E_2$, where
\begin{equation}\label{eq:barApm}
\bar{A}_+ = M^{-1}(-\dot{M}+A_+ M), \quad \hbox{and} \quad \bar{A}_- = M^{-1}(-\dot{M}+A_- M).
\end{equation}
\end{Def}

\begin{Prop}
\label{prop:prod_form}
An optimal synthesis admits a product form if and only if it is the synthesis of a LQ optimal control problem admitting a product structure.
\end{Prop}

\begin{proof}
Note first that, if a problem admits a product structure, then through an appropriate time-dependent linear change of coordinates it splits into $N>1$ sub-problems, and so do the minimizing solutions and the matrices $\bar{A}_+$ and $\bar{A}_-$. This gives the decomposition and proves the if part.

Now, assume that an optimal synthesis characterized by mappings $A_+,A_-$ admits a product form. Up to a (constant) linear change of coordinates, we suppose that the matrices $\bar{A}_+(t)$ and $\bar{A}_-(t)$ have the same block-diagonal form
$\bar{A}_\pm(t)=
\begin{pmatrix}
  A_\pm^1(t) & 0 \\
0 & A_\pm^2(t)
\end{pmatrix}$
for every $t \in \R$, with $A_\pm^i(t) \in M_{n_i}(\R)$, $i=1,2$. Note that, from the definition of $A_\pm$ (given in Theorem~\ref{th:decomp_sol})  and the expression~\eqref{eq:barApm} of $\bar{A}_\pm$, there holds
\begin{equation}
\bar{A}_+ - \bar{A}_- = M^{-1} B R^{-1} B^\top (P_- -P_+)M.
\end{equation}
Since $\bar{A}_+ - \bar{A}_-$ is block diagonal, there exist a smooth map $U: \R \to GL_m(\R)$ such that $\bar{B}=M^{-1} B U$ is also block-diagonal.
As a consequence, through the time-dependent linear change of coordinates
\begin{equation}
\bar{x}=M(t)^{-1}x, \qquad \bar{u}=U(t)^{-1}(u+K(t) x), \quad t\in \R,
\end{equation}
any trajectory of $\dot x = A x + Bu$ is transformed into a trajectory of  $\dot{\bar{x}} = \bar{A}_+ \bar{x} + \bar{B} \bar{u}$. The latter system being block-diagonal, its splits into two independent sub-systems $\dot{\bar{x}}^i = \bar{A}_+^i \bar{x}^i + \bar{B}^i \bar{u}^i$, where $\bar{x}=(\bar{x}^1,\bar{x}^2)$ and $\bar{u}=(\bar{u}^1,\bar{u}^2)$. Moreover the  cost $c(x,u)=(u+Kx)^\top R (u+Kx)$  is equal to $\bar{c}(\bar{x},\bar{u})=\bar{u}^\top \bar{R} \bar{u}$, where $\bar{R}=  U^\top R U$.

Let us decompose the cost as
\begin{equation}
\bar{c}(\bar{x},\bar{u}) = (\bar{u}^1)^\top \bar{R}^1 \bar{u}^1 + (\bar{u}^2)^\top \bar{R}^2 \bar{u}^2 + (\bar{u}^1)^\top \bar{R}^{12} \bar{u}^2.
\end{equation}
If $\bar{R}^{12}=0$, then the system is of the form \eqref{eq:prod} and the proposition is proved. Otherwise, recall that from~\eqref{eq:decomp_sol} and the formula $\Phi_{\bar{A}_\pm}(t,0)=M(t)^{-1} (t) \Phi_{A_\pm}(t,0) M(0)$, every optimal trajectory $\bar{x}(t)$ writes as
\begin{equation}
\bar{x}(t) = \Phi_{\bar{A}_+}(t,t_0) \bar{y}_+ + \Phi_{\bar{A}_-}(t,t_0) \bar{y}_-.
\end{equation}
Since the transition matrices $\Phi_{\bar{A}_\pm}(t,t_0)$ are block-diagonal as well, every optimal trajectory $\bar{x}(t)$ splits into $\bar{x}=(\bar{x}^1,\bar{x}^2)$ where $(\bar{x}^1,0)$ is an optimal trajectory associated with $(\bar{u}^1, 0)$ and $(0,\bar{x}^2)$ is an optimal trajectory associated with $(0,\bar{u}^2)$. It follows that these trajectories minimize also the cost $\bar{c}'(\bar{x},\bar{u}) =(\bar{u}^1)^\top \bar{R}^1 \bar{u}^1 + (\bar{u}^2)^\top \bar{R}^2 \bar{u}^2$, and thus that the optimal synthesis can be obtained by a LQ problem having a product structure.
\end{proof}

\subsection{Injectivity condition}	
\label{se:cns}

\begin{Theorem}
\label{th:noninjectivity}
A canonical cost $c$  admits an equivalent canonical cost if and only if the corresponding LQ optimal control problem admits a product structure.
\end{Theorem}

\begin{proof}
Consider a canonical cost $c=(u+Kx)^\top R (u+Kx)$ and assume that it admits an equivalent canonical cost $\tilde{c}=(u+ \tilde{K} x)^\top \tilde{R} (u+ \tilde{K} x)$.

\begin{Lemma}
	\label{le:injectivity}
If $c$ and $\tilde{c}$ are two equivalent canonical costs associated with $(R,K)$ and $(\tilde{R}, \tilde{K})$ respectively, then $K=\tilde{K}$.
\end{Lemma}

\begin{proof}
If $c \sim \tilde{c}$, then they define the same optimal synthesis $\Gamma$.  From Lemma~\ref{lem:A+A-}, $\Gamma$ determines in a unique way the pair of mappings $(A_+,A_-)$ corresponding to $c$ and $\tilde{c}$, and so a unique pair of mappings $(K_+,K_-)$ since all values of $B$ are injective by \ref{cond:B_injectif}. And Proposition~\ref{prop:reduced} implies that $K_+=K=\tilde{K}$.
\end{proof}

Thus $\tilde{c}=(u+ K x)^\top \tilde{R} (u+ K x)$. The matrix-valued Hamiltonian mapping $\vec{H}$ and $\vec{\tilde{H}}$ associated with $c$ and $\tilde{c}$ respectively  write as
\begin{equation}
\vec{H} =  \begin{pmatrix}
		A_+ & BR^{-1}B^\top\\
		0 & -A_+^\top
		\end{pmatrix}, \qquad
\vec{\tilde{H}}=  \begin{pmatrix}
		A_+ & B \tilde{R}^{-1}B^\top\\
		0 & -A_+^\top
		\end{pmatrix},
\end{equation}
where $A_+=A-BK$.

Lets us prove first a result on the solutions of the corresponding Hamiltonian differential equations.
\begin{Lemma}
\label{prop:orb.diffeo}
There exists a smooth mapping $D: \R \rightarrow GL_n(\R)$ such that $(x, p)$ is a solution of the Hamiltonian differential equation~\eqref{eq:ham_edo} defined by $\vec{H}$ if and only if $(x, Dp)$ is a solution of the Hamiltonian differential equation~\eqref{eq:ham_edo} defined by $\vec{\tilde{H}}$.%
\footnote{Note that the map $(x,p) \mapsto (x,Dp)$ is called an \emph{orbital diffeomorphism} between the Hamiltonians $\vec{H}$ and $\vec{\tilde{H}}$, see~\cite{Jean2019a}.}

Moreover, there holds for every $t \in \R$,
\begin{equation}\label{eq:D(t)}
D(t) = \Phi_{-A_+^\top}(t,0) D(0) \Phi_{-A_+^\top}(t,0)^{-1}.
\end{equation}
\end{Lemma}

\begin{proof}
Let $(x, p)$ be a solution of the Hamiltonian differential equation~\eqref{eq:ham_edo} defined by $\vec{H}$. Then by Theorem~\ref{th:decomp_sol}-(iv), any restriction of $x$ to a subinterval is a minimizer of the LQ problem associated with the cost $c$. Since $c \sim\tilde{c}$, $x$ is also a minimizer  of the LQ problem associated with $\tilde c$, thus there exists $\tilde{p}$ such that $(x,\tilde{p})$ is a solution  of the Hamiltonian differential equation~\eqref{eq:ham_edo} defined by $\vec{\tilde{H}}$. As a consequence,
\begin{equation}
	\dot x = A_+x + BR^{-1}B^\top p = A_+ x + B\tilde{R}^{-1}B^\top \tilde{p}.
\end{equation}
This implies $BR^{-1}B^\top p = B\tilde{R}^{-1}B^\top\tilde{p}$ and, since by Assumption~\ref{cond:B_injectif} all values of $B$ are injective, we obtain
\begin{equation}\label{eq:dotx}
B^\top \tilde{p} = \tilde{R}R^{-1}B^\top p.
\end{equation}

Let us differentiate this equation with respect to time. Using that $\dot p = - A_+^\top p$ and $\dot{\tilde{p}}=-A_+^\top \tilde{p}$, we get
\begin{equation}
\left(\dot{B}^\top - B^\top A_+^\top \right) \tilde{p} = \left( \frac{d}{dt}(\tilde{R}R^{-1}B^\top) - \tilde{R}R^{-1}B^\top A_+^\top \right)  p.
\end{equation}
Set $B_0=B$, $\tilde{B}_0= B R^{-1}\tilde{R}$ and ${B}_{k+1}= \dot{{B}}_k - A_+ {B}_k$, $\tilde{B}_{k+1}= \dot{\tilde{B}}_k - A_+ \tilde{B}_k$ for $k \in \N$.
Taking successive derivatives of the preceding equality, we obtain that, for every $t \in \R$,  $\tilde{p}(t)$ is solution of a system of linear equations,
\begin{equation}
	\label{eq: HtildH}
	\left\{
	\begin{array}{rcl}
	 B_0(t)^\top \tilde{p} (t) & = & \tilde{B}_0 (t)^\top p(t), \\
	& \vdots &\\
	B_k(t)^\top \tilde{p} (t) & = & \tilde{B}_k (t)^\top p(t), \\
& \vdots &\\
\end{array}
	\right.
	\end{equation}
By Assumption~\ref{cond:kalman} and Remark~\ref{re:assumptions}-(3),  this system has a unique solution $\tilde{p}(t)=D(t)p(t)$ for every $t$. Exchanging the role of $c$ and $\tilde{c}$ we obtain that $D(t)$ is moreover invertible, which ends the proof of the first part of the lemma.

From the differential equations $\dot p = - A_+^\top p$ and $\dot{\tilde{p}}=-A_+^\top \tilde{p}$, we get
\begin{equation}
p(t) = \Phi_{-A_+^\top}(t,0) p(0) \quad \hbox{and} \quad \tilde{p}(t) = \Phi_{-A_+^\top}(t,0) \tilde{p}(0).
\end{equation}
Using $\tilde{p}(t)=D(t)p(t)$,   we obtain for every $t$,
\begin{equation}
D(t)\Phi_{-A_+^\top}(t,0)p(0) = \Phi_{-A_+^\top}(t,0) D(0) p(0),
\end{equation}
and then~\eqref{eq:D(t)} since the above equation holds for every $p(0) \in \R^n$.
\end{proof}

Let $D$ be the mapping given by Lemma~\ref{prop:orb.diffeo}. For any trajectory $(x,p)$ of the Hamiltonian differential equation~\eqref{eq:ham_edo} defined by $\vec{H}$, equality~\eqref{eq:dotx} implies $B^\top D p = \tilde{R}R^{-1}B^\top p$. Hence
\begin{equation}\label{eq:D^T}
D^\top B = B R^{-1}\tilde{R}.
\end{equation}
Fix $t \in \R$. Then for every eigenvector $e_t$ of $R^{-1}(t)\tilde{R}(t)$ associated with an eigenvalue $\lambda_t$, the formula above implies that $B(t)e_t$ is an eigenvector of $D(t)^\top$ associated with the eigenvalue $\lambda_t$. Since by~\eqref{eq:D(t)} the eigenvalues of $D(t)$ are the same as the one of $D(0)$, we conclude that the eigenvalues of $R^{-1}(t)\tilde{R}(t)$ are all constant.

As a consequence there exists  a smooth mapping $U: \R \to GL_m(\R)$ such that, for every $t$, $U(t)$ is orthogonal with respect to $R(t)$ and diagonalizes $R(t)^{-1}\tilde{R}(t)$ (see for instance~\cite[Chap.~II, $\S$4]{Kato1995} for the construction of such a mapping $U$). In other terms, there holds $U^\top R U=I$ and $U^{-1}R^{-1}\tilde{R} U= \Lambda$, where $\Lambda$ is a diagonal matrix whose coefficients are the constant eigenvalues of $R^{-1}\tilde{R}$. Note moreover that the hypothesis $c \neq \tilde{c}$ implies that $\Lambda$ admits $N >1$ different eigenvalues $\lambda_1,\dots , \lambda_N$.

Let us apply again Lemma~\ref{prop:orb.diffeo}. Since for any trajectory $(x,p)$ of $\vec{H}$, $(x,Dp)$ is a trajectory of $\vec{\tilde{H}}$, there holds $\frac{d}{dt}(Dp)=-A_+^\top Dp$. Using $\dot p = - A_+^\top p$ we obtain $(\dot{D} - DA_+^\top )p=-A_+^\top Dp$. Since the latter identity holds for every trajectory $(x,p)$, we obtain (the second equality comes from~\eqref{eq:D^T})
	\begin{equation}
	\label{eq:sysonD}
	\dot{D}^\top = A_+D^\top - D^\top A_+ \qquad \hbox{and} \qquad D^\top B U = B U \Lambda.
	\end{equation}
Set  $\bar{B}_0=BU$ and $\bar{B}_{k+1}= \dot{\bar{B}}_k - A_+ \bar{B}_k$ for $k \in \N$. A trivial induction argument shows that $D^\top \bar{B}_k = \bar{B}_k \Lambda$. This implies that, for every $t \in \R$, any column vector of a matrix $\bar{B}_k(t)$, $k \in \N$, is an eigenvector of $\dot{D}(t)^\top$ associated with an eigenvalue equal to one of the $\lambda_1,\dots , \lambda_N$. Since by Assumption~\ref{cond:kalman} and Remark~\ref{re:assumptions}-(3) these column vectors generate $\R^n$, we deduce that ${D}(t)^\top$ is diagonalizable for every $t$, with spectrum equal to $\{ \lambda_1,\dots , \lambda_N \}$.

This diagonalization can actually be realized in a basis depending smoothly on $t$. Indeed, let $M_0 \in GL_n(\R)$ be a matrix diagonalizing $D(0)^\top$ and set $M(t)=\Phi_{-A_+^\top}(t,0)^{-\top} M_0$, $t \in \R$. Then $M: \R \to GL_n(\R)$ is a smooth mapping and, by~\eqref{eq:D(t)}, for every $t \in \R$ we have
\begin{equation}
  M(t)^{-1} D(t)^\top M(t) = \Delta,
\end{equation}
where (up to reordering) $\Delta$ is a diagonal matrix of the form
\begin{equation}
\Delta = \begin{pmatrix}
	\lambda_1 I_{n_1} & & \\
	& \ddots & \\
	& & \lambda_N I_{n_N}
	\end{pmatrix}.
\end{equation}
\smallskip

Let us now use the reductions above to construct a change of state and control coordinates putting the optimal control problem into product form. Consider the time-dependent linear transformation $G_t : (x,u) \in \R^n \times \R^m \to (\bar{x},\bar{u}) \in \R^n \times \R^m$ defined by
\begin{equation}
\bar{x}=M(t)^{-1}x, \qquad \bar{u}=U(t)^{-1}(u+K(t) x), \quad t\in \R.
\end{equation}
The image $(\bar{x}(t),\bar{u}(t))=G_t(x(t),u(t))$, $t \in \R$, of any trajectory of $\dot x = A x + Bu$ satisfies  $\dot{\bar{x}} = \bar{A} \bar{x} + \bar{B} \bar{u}$ where
\begin{equation}
\bar{A} = M^{-1}(-\dot{M}+A_+ M) \quad \hbox{and} \quad \bar{B} = M^{-1}BU,
\end{equation}
and its quadratic cost writes as
\begin{equation}
  c(t,x(t),u(t)) = \bar{u}(t)^\top \bar{u}(t).
\end{equation}
We shall have established the theorem if we prove that $\bar{A}$, $\bar{B}$ are of the form~\eqref{eq:prod}.

Let $E_1, \dots, E_N$ be the eigenspaces of $\Delta$. Let us remark that
\begin{equation}
  0= \dot{\Delta}=\frac{d}{dt}(M^{-1} D^\top M) = \bar{A} \Delta - \Delta \bar{A}.
\end{equation}
Thus for every $t \in \R$, the matrices $\bar{A}(t)$ and $\Delta$ commute, which implies that $\bar{A}(t)$ preserves every $E_j$. On the other hand, by~\eqref{eq:sysonD} every column of $B(t)U(t)$ is an eigenvector of $D(t)^\top$, which implies that every column of $\bar{B}(t) = M(t)^{-1}B(t)U(t)$ is an eigenvector of $\Delta$, and then belongs to one of the eigenspaces $E_j$. Thus, in the decomposition $\R^n = E_1 \oplus \dots \oplus E_N$, the matrices $A$ and $B$  have the required product form,
	\begin{equation*}
	\bar{A} =
	\begin{pmatrix}
	A_1 & & \\
	& \ddots & \\
	& & A_N
	\end{pmatrix}
	\quad
	B =
	\begin{pmatrix}
	B_1 &  & \\
	& \ddots & \\
	&  & B_N  \\
	\end{pmatrix},
	\end{equation*}
which ends the proof.
\end{proof}

As a consequence of this theorem and of Proposition~\ref{prop:prod_form}, we obtain the following characterization of the uniqueness of solutions to the inverse LQ problem.
\begin{Coro}
\label{cor:m>1}
An inverse optimal control problem among canonical costs admits a unique solution if and only if the optimal synthesis does not admit a product form.
\end{Coro}

Since the number $N$ of elements in a product structure satisfies $1<N\leq m$, we obtain in particular an injectivity result in the single input case (this generalizes the result of \cite{Nori2004} to non autonomous problems).
\begin{Coro}
\label{cor:m=1}
In the single input case ($m=1$), any inverse optimal control problem in canonical form admits a unique solution.
\end{Coro}


\subsection{Reconstruction}
\label{se:reconstruction_nonaut}
Let us consider now the problem of the reconstruction of the cost. We fix a pair $(A,B)$ satisfying Assumption~\ref{Assump:AB}, and we look for an explicit solution to the following inverse problem: given an optimal synthesis, recover a canonical cost that would produce this optimal synthesis. Such a construction includes two steps:
\begin{itemize}
  \item identify the unique pair of mappings $(A_+,A_-)$ associated with the optimal synthesis;
  \item construct the mappings $(R,K)$ defining a canonical cost solution to the problem  as functions of $(A_+,A_-)$.
\end{itemize}

The first step is essentially a numerical identification problem, it will be treated in the autonomous case in Section~\ref{se:numerics}. Let us focus here on the second step.

\begin{Prop} \label{prop:reconstr_nonautonome}
Let $(A_+,A_-)$ be a pair of mapping describing an optimal synthesis. The optimal costs  whose optimal synthesis is described by $(A_+,A_-)$ are the costs $c=(u+Kx)^\top R (u+Kx)$ such that:
\begin{enumerate}
  \item the mapping $K$ is given by
\begin{equation}
  K = \left(B^{\top}B\right)^{-1}B^{\top}\left(A - A_+ \right) ;
\end{equation}
  \item the mapping $R$ is given by $R(t)=(\det Z(0))^{1/n} Z^{-1}(t)$, where $Z: \R \to \Sym (m)(\R)$ is a smooth mapping with positive definite values solution of the linear equation
\begin{multline}\label{eq:Z_nonautonome}
 (A_- - A_+) \int_{0}^{t} \Phi_{A_+}(s,0) B(s) Z(s) B(s)^\top \Phi_{A_+}(s,0)^\top ds \\
 - B(t) Z(t) B(t)^\top (I - \Phi_{A_-}(t,0)^{-\top}\Phi_{A_+}(t,0))=0
\end{multline}
\end{enumerate}
\end{Prop}

\begin{rk}
It results from Corollary~\ref{cor:m>1} that, if the optimal synthesis does not admit a product form, then the linear equation~\eqref{eq:Z_nonautonome} admits a unique solution $Z$ such that $\det Z(0)=1$. In that case, the above proposition allows to construct the unique canonical cost solution to the inverse problem by solving a linear equation.
\end{rk}

\begin{proof}
Let $c=(u+Kx)^\top R (u+Kx)$ be a canonical cost whose optimal synthesis is described by $(A_+,A_-)$. By Proposition~\ref{prop:reduced} there holds $A_+=A-BK$, which directly implies the expression of $K$ (recall that all values of $B$ are injective by \ref{cond:B_injectif}).

Now, from Theorem~\ref{th:decomp_sol} there exists mappings $P_+,P_-$ such that $A_\pm = A - BR^{-1} (RK + B^\top P_\pm)$ and $P_+ - P_-$ has positive definite values. Setting $X= ( P_+ - P_-)^{-1}$ and $Z=R^{-1}$, we obtain
\begin{equation}\label{eq:A-moinsA+}
 (A_- - A_+) X = B Z B^\top.
\end{equation}

On the other hand, $P_+,P_-$ are also solution of the associated Riccati differential equation. A direct (and rather standard, see~\cite[Sect.\ 4.1]{AbouKandil2003}) computation shows that $X$ is solution of the following Lyapunov differential equations,
\begin{equation}\label{eq:X}
 \dot X = A_+X + X A_+^\top + B Z B^\top \quad \hbox{and} \quad \dot X = A_+X + X A_-^\top.
\end{equation}
The solution of these differential equations admit explicit expressions, respectively
\begin{equation}
 X (t) = \Phi_{A_+}(t,0) X(0) \Phi_{A_+}(t,0)^\top + \int_{0}^{t} \Phi_{A_+}(s,0) B(s) Z(s) B(s)^\top \Phi_{A_+}(s,0)^\top ds,
\end{equation}
\begin{equation}
\hbox{and} \quad X (t) = \Phi_{A_+}(t,0) X(0) \Phi_{A_-}(t,0)^\top .
\end{equation}
By plugging in these expressions into~\eqref{eq:A-moinsA+} we obtain two different equations. Eliminating $X(0)$ between these equations, we obtain~\eqref{eq:Z_nonautonome}.
\end{proof}

\section{Autonomous case}
\label{sec:autonomous}
Let us now consider  in this section the \emph{autonomous case}, that is the particular case where $A,B,Q,S,R$ are constant matrices. All results obtained in the previous sections may then be reformulated in an autonomous setting (with only constant matrices), and under weaker assumptions on $(A,B,Q,S,R)$.

\begin{assumption}
	\label{Assump_aut:AB}
	The matrices $(A,B)$ satisfy:
	\begin{enumerate}[label=\textit{(\ref{Assump_aut:AB}.\arabic*)}]
		\item  \label{cond:C1} $B$ is of rank $m$;

		\item  \label{cond:C2}   $(A,B)$ satisfy the Kalman controllability condition
		\begin{equation}
			\mathrm{span}\{B, AB, \dots, A^{n-1}B\} = \R^n.
		\end{equation}
	\end{enumerate}
\end{assumption}

\begin{assumption}
	\label{Assump_aut:ABQSR}
	The matrices $(Q,S,R)$ satisfy:
	\begin{enumerate}[label=\textit{(\ref{Assump_aut:ABQSR}.\arabic*)}]
		\item  \label{cond:D1} $Q \in \Sym(n), \ R \in \Sym(m)$ and  $R \succ 0$; 
		
		\item  \label{cond:D2}  the Hamiltonian matrix
		$$
		\begin{pmatrix}
			A - BR^{-1}S^\top & BR^{-1}B^\top\\
			Q - SR^{-1}S^\top & -A^\top + SR^{-1}B^\top
		\end{pmatrix}
		$$
		has no eigenvalues on the imaginary axis.
	\end{enumerate}
\end{assumption}
In this setting, the optimal synthesis may be characterized through the solutions of the algebraic Riccati equation,
\begin{equation}
	\label{eq_alg:Riccati}
	PA + A^\top P - (S +PB) R^{-1}(S^\top + B^\top P) + Q = 0,
\end{equation}
rather than through the ones of the differential Riccati equation, and the characterization of Theorem~\ref{th:decomp_sol} is replaced by the following one.

\begin{Theorem}
	\label{th:decomp_sol_aut}
	Consider matrices $(A,B,Q,S,R)$ satisfying Assumptions~\ref{Assump_aut:AB} and~\ref{Assump_aut:ABQSR}.
	
	\begin{enumerate}[label={(\roman*)}]
		\item \label{th:(1)} There exists a unique positive definite $P_+ \in \Sym(n)$ maximal solution of the algebraic Riccati equation~\eqref{eq_alg:Riccati} such that the matrix
		\begin{equation}
			\label{eq:A+.aut}
			A_{+} = A-B R^{-1} \left(S^\top + B^\top P_{+}\right)
		\end{equation}
		is asymptotically stable, i.e. has only eigenvalues with negative real part.
		
		\item \label{th:(2)} There exists a unique negative definite $P_- \in \Sym(n)$ minimal solution of the algebraic Riccati equation~\eqref{eq_alg:Riccati} such that the matrix
		\begin{equation}
			\label{eq:A-.aut}
			A_{-} = A-B R^{-1} \left(S^\top + B^\top P_{-}\right)
		\end{equation}
		is asymptotically anti-stable, i.e. has only eigenvalues with positive real part.

		\item \label{th:(3)} For every data $t_0<t_1$, $x_0,x_1 \in \R^n$, there exists a unique trajectory $x(\cdot)$ of the corresponding autonomous LQ problem associated with $(A,B,Q,S,R)$. Moreover there exists a unique pair $(y_+,y_-) \in \R^{2n}$ such that
		\begin{equation}
			\label{eq:decomp_sol_aut}
			x(t) =  e^{tA_{+}} y_+ + e^{tA_{-}} y_-, \qquad \hbox{for every } t \in [t_0,t_1].
		\end{equation}
	\end{enumerate}
\end{Theorem}

\begin{proof}
	Parts \textit{(i)} and \textit{(ii)} result from \cite[Th.\ 6]{Molinari1977}. The existence of minimizers  in \textit{(iii)} is a consequence of \cite[Th.\ A]{Agrachev2015}, and \eqref{eq:decomp_sol_aut} is obtained in the same way as~\eqref{eq:decomp_sol}, see section~\ref{se:appendix}.
\end{proof}

\begin{rk}
	\begin{itemize}
		\item Assumption~\ref{cond:C1} is actually not necessary for this result to hold. It will be necessary however to address the inverse problem.
		\item Assumptions~\ref{cond:C2} and~\ref{cond:D1} are necessary for existence and uniqueness of the minimal solution, but not sufficient as stated in~\cite[Th.\ A]{Agrachev2015}. On the other hand, Assumption~\ref{cond:D2} is necessary and sufficient for the existence of $P_+, P_-$.
		\item This result generalizes \cite[Th.\ 1]{Ferrante2005}, which is obtained under the extra assumption that $\bar{Q}= Q - S R^{-1} S^\top$ is positive semi-definite.
	\end{itemize}
\end{rk}

\subsection{Inverse autonomous LQ optimal control problem}
\label{se:injectivity_aut}
The inverse autonomous linear-quadratic optimal control problem is posed as follows: a pair of matrices $(A,B)$ satisfying Assumption~\ref{Assump_aut:AB} being fixed, recover the autonomous quadratic cost $c$ from the given of the corresponding optimal synthesis. So the only difference with respect to the previous sections is that we are looking for solutions of the inverse problem in the smaller class of autonomous costs. It is straightforward to verify that all preceding results still hold with ``autonomous'' statements (and actually most of the results were already obtained in \cite{Jean2018a} under the extra hypothesis that $\bar{Q}= Q - S R^{-1} S^\top$ is positive semi-definite). We therefore state these results without proof, simply highlighting when needed the difference with the non-autonomous case.

The injectivity problem induces the equivalence relation on costs $(Q,S,R)$ as in Definition~\ref{def:equiv}. The analog of Lemma~\ref{lem:A+A-} in the autonomous case allows us to state the equivalence relation on the pairs of matrices $(A_+, A_-)$ defined by \eqref{eq:A+.aut} and \eqref{eq:A-.aut} in Theorem~\ref{th:decomp_sol_aut}.

\begin{Lemma} \label{Lemma:A+A-}
	Two equivalent autonomous costs define the same matrices $(A_+, A_-)$.
\end{Lemma}
Equivalently, two equivalent autonomous costs define the same matrices $K_+, K_- \in M_{n,m}(\R)$, where $A_{+} = A-BK_{+}$ and $A_{-} = A-BK_{-}$.

Canonical autonomous costs are defined  as in Definition~\ref{def:KR} by
\begin{equation} \label{eq:canonical.cost.aut}
	c(x,u)= (u+K x)^\top R (u+Kx),
\end{equation}
where
\begin{itemize}
	\item $R \in \Sym(m)$ is a positive definite matrix such that  $\det R=1$;
	\item the matrix $K \in M_{m,n}(\R)$ is such that $A - BK$ is asymptotically stable.
\end{itemize}
Any autonomous cost $(Q,S,R)$ has an equivalent cost in canonical form with $K = K_+$, therefore the inverse problem can  be stated on the class of canonical costs. The adaptation of Definition~\ref{def:product} to the autonomous case is as follows.
\begin{Def}
	We say that an autonomous LQ optimal control problem defined by a pair of matrices $(A,B)$ and an autonomous canonical cost $c=(u+Kx)^\top R (u+Kx)$ \emph{admits a product structure} if there exists an integer $N> 1$ and matrices  $M \in GL_n(\R)$ and $U \in GL_m(\R)$ such that
	\begin{equation}
		\bar{A} = M^{-1}A M, \quad \bar{B} = M^{-1}BU, \quad \hbox{and} \quad \bar{R} = U^\top R U
	\end{equation}
	are of the form \eqref{eq:prod}.
\end{Def}
The non injective  inverse autonomous optimal problems are the ones corresponding to LQ problems admitting a product structure, as stated in this counterpart of Theorem~\ref{th:noninjectivity}.

\begin{Theorem}
	\label{th:prod_aut}
	An inverse autonomous LQ optimal control problem among canonical costs admits a unique solution if and only if the synthesis is not obtained by a problem having a product structure.
\end{Theorem}

\subsection{Reconstruction in the autonomous case}
\label{se:reconstruction}
Let us fix a pair of matrices $(A,B)$ satisfying Assumption \ref{Assump_aut:AB}. We  look for an explicit solution to the following inverse problem: given an optimal synthesis, recover a canonical cost that would produce this optimal synthesis. As mentioned in Section~\ref{se:reconstruction_nonaut} such a construction includes two steps:
\begin{enumerate}[label={(\roman*)}]
	\item identify the unique pair of matrices $(A_+,A_-)$ associated with the optimal synthesis;
	\item construct the matrices $(R,K)$ defining a canonical cost solution to the problem  as functions of $(A_+,A_-)$.
\end{enumerate}
Let $\mathcal{J}$ be the set of pairs $(A_+, A_-)$ in $ M_n(\R)^2$ corresponding to optimal synthesis of autonomous LQ problem. To identify algorithmically $(A_+,A_-)$ from the set of trajectories, we need first to characterize $\mathcal{J}$. 

Consider an autonomous LQ problem and set $\Delta = P_+ - P_-$, where $P_+,P_-$ are the matrices given by Theorem~\ref{th:decomp_sol_aut}. Since $\Delta$ is a positive definite matrix, it inverse  $X = \Delta^{-1}$ is well defined. It follows from \eqref{eq_alg:Riccati} that
$X$ satisfies the following algebraic Lyapunov equations,
\begin{equation} \label{eq:lyap_alg}
	A_+ X + XA_+^\top = - B R^{-1} B^\top \quad \hbox{ and } \quad A_+ X + XA_-^\top = 0,
\end{equation}
the solution to the left one being given by
\begin{equation} \label{eq:solX}
	X = \int_{0}^{\infty}e^{t A_+} BR^{-1}B^{\top} e^{t(A_+)^{\top}} dt.
\end{equation}

\begin{Prop}
	The set $\mathcal{J}$ is equal to the set of pairs $(A_+ = A-BK_+, A_- = A-BK_-)$, with $K_\pm \in M_{m,n}(\R)$, for which there exists $\Delta \in \Sym^+(n)$ such that:
	\begin{itemize}
		\item $A_+ = - \Delta^{-1} A_-^\top  \Delta$;
		\item $A_+$ is asymptotically stable;
		\item $(B^\top B)^{-1}B^\top\left( A_+ \Delta^{-1} + \Delta^{-1}A_+^\top \right)B(B^\top B)^{-1} \in  \Sym^-(m)$.
	\end{itemize}
\end{Prop}


%
\begin{proof}
	Denote by $\mathcal{A}$ the set of pairs satisfying the properties above.  By definitions \eqref{eq:A+.aut}-\eqref{eq:A-.aut} and by  \eqref{eq:lyap_alg}, any $(A_+, A_-)$ in $\mathcal{J}$ belongs to $\mathcal{A}$ as well. Conversely, for any $(A_+, A_-)$ in $\mathcal{A}$, let us define $K$ as the unique solution of  $A_+=A-BK$, and $R$ by
	\begin{equation} \label{eq:reconstruction.R}
		R^{-1} = - (B^\top B)^{-1}B^\top\left( A_+ \Delta^{-1} + \Delta^{-1}A_+^\top \right)B(B^\top B)^{-1}.
	\end{equation}
	Then $(A_+, A_-)$  coincide with the matrices \eqref{eq:A+.aut},\eqref{eq:A-.aut} associated with the LQ problem with canonical cost determined by $K$ and $R$, which proves that $\mathcal{A} \subset \mathcal{J}$.
\end{proof}

The inverse problem  is not injective on the pairs $(A_+, A_-)$ corresponding to LQ problems with product structure. Such pairs admit the following simple characterization, which is an adaptation of Proposition~\ref{prop:prod_form}.

\begin{Prop}
	\label{prop:prodA+A-}
	An autonomous LQ optimal control problem admits a product structure if and only if the associated matrices $A_+$ and $A_-$  both leave invariant  a nontrivial (i.e.\ $E_1,E_2 \neq \{0\}$) decomposition $\R^n=E_1 \oplus E_2$.
\end{Prop}

Let us denote by $\mathcal{J}_p$ the subset of $\mathcal{J}$ composed by the pairs $(A_+, A_-)$ in $\mathcal{J}$ which leave invariant a nontrivial decomposition of $\R^n$. By Theorem~\ref{th:prod_aut} and Proposition~\ref{prop:prodA+A-}, the inverse optimal control problem in canonical form is injective on $\mathcal{J} \setminus \mathcal{J}_p$. We will prove in next section that $\mathcal{J} \setminus \mathcal{J}_p$ contains an open and dense set in $\mathcal{J}$, which means that a generic optimal synthesis corresponds to a unique canonical cost. In practice, with noisy data, one can always assume that we are in the situation of generic optimal synthesis, and thus, it is enough to treat the reconstruction for the injective case only. This reconstruction is given in  Algorithm~\ref{alg:reconstruction}. 

\begin{algorithm}
	\caption{Reconstruction of $(R, K)$ from an optimal synthesis $\Gamma$}\label{alg:reconstruction}
	\begin{algorithmic}
		\State 1. Identify from the trajectories in $\Gamma$ the matrices $(K_+, K_-, \Delta)$ characterizing $(A_+, A_-)$ in $\mathcal{J}$.
			%
			\State 2. Set $K = K_+$
			\State 3. Compute $Z=-(B^\top B)^{-1}B^\top\left( A_+ \Delta^{-1} + \Delta^{-1}A_+^\top \right)B(B^\top B)^{-1}$.		
			\State 4. Set  $R=(\det Z)^{1/n} Z^{-1}$.
		\end{algorithmic}
	\end{algorithm}

\subsection{Genericity of the injective case}
\label{sse:genericity}

Fix matrices $A$, $B$ satisfying Assumptions~\ref{Assump_aut:AB}. Recall that $\mathcal{J}$ is the set of pairs $(A_+, A_-)$ in $ M_n(\R)^2$ corresponding to optimal synthesis of autonomous LQ problem, and $\mathcal{J}_p$ is the subset of $\mathcal{J}$ composed of $(A_+, A_-)$ which leave invariant a nontrivial decomposition of $\R^n$. We equip $\mathcal{J}$ with the induced topology of $M_n(\R)^2$. 

\begin{Theorem} \label{th:genericity}
	The set $\mathcal{J} \setminus \mathcal{J}_p$ contains an open and dense subset of  $\mathcal{J}$.
\end{Theorem}
\begin{proof}
	Let $\Omega$ be the set of pairs $(A_+, A_-)$ in $M_n(\R)^2$ such that $A_+$ and $A_-$ have only simple eigenvalues and do not leave invariant any nontrivial decomposition of $\R^n$. By definition, the intersection of $\mathcal{J} \cap \Omega$ is included in $\mathcal{J} \setminus \mathcal{J}_p$.  We will show that this intersection is an open and dense subset of  $\mathcal{J}$.
	
	Let us first show  that $\Omega$ is open in $M_n(\R)^2$, which will imply that $\mathcal{J} \cap \Omega$ is open in $\mathcal{J}$.  Consider the subset $\Omega_s$ of $ M_n(\R)^2$ consisting of pair of matrices having both only simple eigenvalues. This open set $\Omega_s$ contains $\Omega$. Let $(A_+, A_-)$ be an element of $\Omega_s$, and let $e_1^\pm, \dots, e_n^\pm$ be a basis of $\C^n$ made of eigenvectors of $A_\pm$ (both matrices $A_\pm$ are diagonalizable since their eigenvalues are simple). The matrices $A_+, A_-$ leave invariant a nontrivial decomposition $E_1 \oplus E_2 = \R^n$ if and only if there exists two nontrivial partitions $I_+ \cup J_+ = \{1, \dots, n\}$ and $I_- \cup J_- = \{1, \dots, n\}$ such that
	\begin{equation}
		\left\{
		\begin{array}{l}
			E_1 = \R^n \cap \mathrm{span}\{e_i^+ \ : i \in I_+\} = \R^n \cap \mathrm{span}\{e_i^- \ : i \in I_-\}   , \\[2mm]
			E_2 = \R^n \cap \mathrm{span}\{e_j^+ \ : j \in J_+\} = \R^n \cap \mathrm{span}\{e_j^- \ : j \in J_-\}.
		\end{array}
		\right.
	\end{equation}
	This implies that $\langle e^+_i , e^-_j \rangle_{A_+} = 0$ for all $(i,j)$ in $I_+ \times J_-$ or in $I_- \times J_+$, where $\langle \cdot, \cdot \rangle_{A_+}$ denotes the Hermitian product associated with $(e^+_1, \dots, e^+_n)$. Therefore, the pair $(A_+, A_-)$ belongs to $\Omega$ if and only if for any nontrivial partitions $I_\pm \cup J_\pm$ of $\{1, \dots, n\}$ there exists $(i,j)$ in  $I_+ \cup J_-$ or in $J_+ \cup I_-$ such that $\langle e^+_i , e^-_j \rangle_{A_+} \neq 0$. Now, the eigenvectors $e^\pm_i$, and then the corresponding Hermitian product $\langle \cdot, \cdot \rangle_{A_+}$, can be chosen locally as continuous functions of $(A_+, A_-)$ in $\Omega_s$. Thus the condition for a pair $(A_+, A_-)$ in $\Omega_s$ to belong to $\Omega$ is an open condition. This implies that $\Omega$ is an open subset of $\Omega_s$, and thus an open subset of $M_n(\R)^2$.\smallskip
	
	Let us show now that the intersection $\mathcal{J} \cap \Omega$ is dense in $\mathcal{J}$. We begin by writing $\mathcal{J}$ as the image of a function. Let $\Sym^+(m)$ be the set of positive definite matrices in $\Sym(m)$ and define the continuous function $\Phi : M_{n,m}(\R) \times \Sym^+(m) \to M_n(\R)^2$ by $\Phi(K,R)= (A_+, A_-)$, where
	\begin{equation}
		A_+ = A - BK, \quad A_- = -X^{-1} A_+ X, \quad \hbox{where } X = \int_{0}^{\infty}e^{t A_+} B R^{-1}B^{\top} e^{t(A_+)^{\top}} dt.
	\end{equation}
	Notice that $\mathrm{sp}(A_-) =- \mathrm{sp}(A_+)$ by construction.
	By \eqref{eq:lyap_alg} and~\eqref{eq:solX}, we have $\mathcal{J} = \Phi(\mathcal{M})$ where $\mathcal{M}$ is the open subset of $M_{n,m}(\R) \times \Sym^+(m)$ defined by
	$$
	\mathcal{M} = \{ (K, R) \in  M_{n,m}(\R) \times \Sym^+(m)  \ : \ A - BK \text{ is asymptotically stable} \}.
	$$
	Let us denote by $\mathcal{M}_s$ the subset of $\mathcal{M}$ containing the pairs $(K, R)$ such that $A - BK$ has simple spectrum and the LQ problem associated with $A-BK$, $B$, and $R$ does not have a product structure (i.e.\ the matrices can not be simultaneously block-diagonalized by a change of coordinates of $\R^n$ and $\R^m$). Since the matrices $A_+$ and $A_-$ defining $\Phi(K, R)$ have opposite spectrum, by Proposition~\ref{prop:prodA+A-} there holds
	\begin{equation}
		\Phi(\mathcal{M}_s) = \mathcal{J} \cap \Omega.
	\end{equation}
	The map $\Phi$ being continuous, it suffices to show that $\mathcal{M}_s$ is dense in $\mathcal{M}$ in order to show that $\mathcal{J} \cap \Omega$ is dense in $\mathcal{J}$.
	
	Let $(K^0, R^0)$ be an element in $\mathcal{M}$ and $\mathcal{W} \subset \mathcal{M}$ be a neighbourhood of $(K^0, R^0)$. We claim that:
	
	\begin{enumerate}
		\item there exists $K^1 \in M_{n,m}(\R)$ such that $(K^1, R^0)$ belongs to $\mathcal{W}$ and $A-BK^1$ has a simple spectrum; consequently, $(K^1, R)$ belongs to $\mathcal{W}$ for every $R$ in an open neighbourhood $\mathcal{N}$ of $R^0$;
		\item there exists $R^1 \in \mathcal{N}$ such that the matrices $A-BK^1$, $B$, and $R^1$ can not be simultaneously block-diagonalized.
	\end{enumerate}
	As a result, we would obtain that the pair $(K^1, R^1)$ belongs to $\mathcal{W} \cap \mathcal{M}_s$, which proves that $\mathcal{M}_s$ is dense in $\mathcal{M}$ and hence ends the proof of the theorem.\smallskip
	
	It remains to prove the claim. Denote by $\mathcal{P}_n$ the vector space of unitary polynomials of degree $n$ with real coefficients, and consider the application $\psi: M_{n,m}(\R) \to \mathcal{P}_n$ which maps a matrix $K$ to the characteristic polynomial of $A-BK$. By the pole shifting theorem (see for instance \cite[Th.~10.1]{Coron2007}), which holds by Assumption \eqref{Assump_aut:AB}, the map $\psi$ is surjective. Actually it is shown in the proof of the latter theorem that the differential of $\psi$ at a certain $K$ is surjective (see bottom of page 277 in \cite{Coron2007}). Since $\psi$ is an algebraic map, we deduce that its differential is surjective at almost all $K \in M_{n,m}(\R)$. This implies that we can find a matrix $K$ arbitrarily close from $K^0$ such that $\psi$ is a locally open map at $K$. As the set of polynomials with simple roots is dense in $\mathcal{P}_n$, we conclude that there exists $K^1$ arbitrarily close from $K^0$ such that the characteristic polynomial of $A-BK^1$ has simple roots. This proves (1).
	
	Denote by $\R^n=E_i \oplus E_i'$, $i=1,\dots, k_n$ all the decompositions of $\R^n$ into vector subspaces invariant by $A-BK^1$ (the number $k_n$ of such decompositions is finite), and denote by $\{1,\dots,m\} = J_j \cup J_j'$, $j=1,\dots,k_m$, all partitions of $\{1,\dots,m\}$ into complementary nonempty subsets.  Given $R \in \Sym^+(m)$, the matrices $A-BK^1$, $B$, and $R$ can be simultaneously block-diagonalized if and only if there exist $i \leq k_n$, $j \leq k_m$, and $U \in GL_m(\R)$ such that
	\begin{itemize}
		\item if $\ell$ in $J_j$ (resp.\ in $J'_j$), then the $\ell$-th column of $BU$ belongs to $E_i$ (resp.\ $E'_i$);
		\item  if $(k,\ell) \in J_j \times J'_j$, then the $k,\ell$ component $(U^\top RU)_{k,\ell}$ equals $0$.
	\end{itemize}
	The pair of indices $(i,j)$ characterizes the form of the resulting block-diagonal matrices $M^{-1}(A-BK^1)M$, $M^{-1}BU$, and $U^\top RU$, where the change-of-basis $M$ can be any matrix whose column are formed by the union of bases of $E_i$ and $E'_i$. There is a finite number $N(R)$ of such diagonalizing pairs of indices at each $R$. Moreover, by a continuity argument, the diagonalizing pairs of indices at $R$ are the only such pairs at $R'$ close enough to $R$.
	
	Let $R^1$ be a matrix realizing the minimal value of $N(R)$ on $\mathcal{N}$ and assume by contradiction that $N(R^1)>0$, i.e., that $A-BK^1$, $B$, and $R^1$ can be simultaneously block-diagonalized. There exists a block-diagonal form defined by a pair $(i_0,j_0)$ such that, up to change of basis in $\R^n$ and $\R^m$,
	\begin{equation}\label{eq:blockdiagonalform}
		A-BK^1 =
		\begin{pmatrix}
			A_1 & 0 \\
			0 & A_2
		\end{pmatrix}, \quad B = \begin{pmatrix}
			B_1 & 0 \\
			0 & B_2
		\end{pmatrix}, \quad R^1 = \begin{pmatrix}
			R_1 & 0 \\
			0 & R_2
		\end{pmatrix}.
	\end{equation}
	Consider a matrix $R^\eps \in \Sym^+(m)$ of the form
	\begin{equation}
		R^\eps = \begin{pmatrix}
			R_1 & \eps R_3 \\
			\eps R_3^\top & R_2
		\end{pmatrix},
	\end{equation}
	where $R_3$ is a nonzero matrix of appropriate dimension and $\eps$ a real number. For $\eps>0$ small enough, $R^\eps \in \mathcal{N}$, $N(R^\eps)=N(R^1)$, and the diagonalizing pairs of indices at $R^\eps$ are necessarily the same as the ones at $R^1$. It is in particular the case for the pair $(i_0,j_0)$, which implies that there exists a matrix $U \in GL_m(\R)$ such that $BU$ and $U^\top R^\eps U$ have the same block-diagonal form~\eqref{eq:blockdiagonalform} as $B$ and $U$. Writing $BU$ as
	\begin{equation}
		BU = \begin{pmatrix}
			B_1 U_{11} & B_1 U_{12} \\
			B_2 U_{21} & B_2 U_{22}
		\end{pmatrix}, \quad \hbox{where } U = \begin{pmatrix}
			U_{11} &  U_{12} \\
			U_{21} & U_{22}
		\end{pmatrix},
	\end{equation}
	we obtain $B_1 U_{12}=0$ and $B_2 U_{21}=0$, which gives $U_{12}=0$ and $U_{21}=0$ (recall that $B_1$ and $B_2$ are injective by Assumption~\ref{cond:C1}). Moreover, $U^\top R^\eps U$ writes as
	\begin{equation}
		U^\top R^\eps U = \begin{pmatrix}
			U_{11}^\top R_1 U_{11} & \eps U_{11}^\top R_3 U_{22} \\
			\eps U_{22}^\top R_3^\top U_{11} & U_{22}^\top R_2 U_{22}
		\end{pmatrix} \quad \hbox{since } U = \begin{pmatrix}
			U_{11} &  0 \\
			0 & U_{22}
		\end{pmatrix},
	\end{equation}
	and has the same block-diagonal form as $R^1$, which implies $U_{11}^\top R_3 U_{22}=0$. This contradicts $U \in GL_m(\R)$ and $R_3 \neq 0$. Thus $N(R^1)=0$, which means that the matrices $A-BK^1$, $B$, and $R^1$ can not be simultaneously block-diagonalized. This shows point (2) of the claim and ends the proof.
\end{proof}

\section{Numerical experiments}
\label{se:numerics}

In this section, we apply Algorithm~\ref{alg:reconstruction} to an example of autonomous LQ problem in $\R^3$. The goal in the considered example is to reconstruct the matrices $K, R$ from the given set of trajectories $\Gamma$. To generate $\Gamma$ for the reconstruction, the direct linear quadratic problems is solved numerically using Matlab and \eqref{eq:decomp_sol_aut}, \eqref{eq:A+.aut}, \eqref{eq:A-.aut}, where $P_+$ is obtained using Matlab function care.m and $\Delta = P_+ - P_-$ using lyap.m. We get the optimal trajectories for 3 variants of boundary conditions, namely, the initial condition given by the vector $x(0) = (0,0,0)$ and 3 variants of final conditions $x_i(T)=e_i$ where $e_i$ denotes the $i$th vector of the canonical basis of $\R^3$. The final time is fixed to $T=1$. By Lemma~\ref{lem:A+A-}, the chosen trajectories determine uniquely a pair $(A_+, A_-)$, and therefore, are enough for the reconstruction.

In practice, the reconstruction is done from a finite number of points of trajectories in $\Gamma$. In our tests we use the points at time grid $t_i$ with steps size $t_{i+1} - t_i = T/N, \ i = 1, \dots, N$ and $N = 20$. In the numerical example, we also test the stability of Algorithm~\ref{alg:reconstruction} with respect to the noise.  The noise is given by $\alpha \times \mathrm{randn}$, where randn.m is a Matlab function generating a normally distributed random values and $\alpha$ is the noise amplitude which ranges from 0 to 0.2 in our tests. As a result of the reconstruction, we obtain a pair $(K_{\mathrm{est}},R_{\mathrm{est}})$ from given trajectories $\Gamma$, which determine a canonical cost \eqref{eq:canonical.cost.aut}.

To evaluate the results, we compare the estimated $(K_{\mathrm{est}},R_{\mathrm{est}})$ with $(K_+, R)$ of the matrices in LQ problem used for the generation of $\Gamma$ in case of 100 samples of noisy data for each value of the amplitude $\alpha = 0, \ 0.05,\ 0.1,\ 0.15,\ 0.2$. The comparison is done by calculating the mean value of $(K_{\mathrm{est}},R_{\mathrm{est}})$ over $L_{\alpha}$ successfully obtained results of the same amplitude $\alpha$
$$
\bar{K}_{\mathrm{est}} = \frac{1}{L_{\alpha}}\sum_{i=1}^{L_{\alpha}}K^i_{\mathrm{est}} \qquad \bar{R}_{\mathrm{est}} = \frac{1}{L_{\alpha}}\sum_{i=1}^{L_{\alpha}}R^i_{\mathrm{est}}.
$$
And then we calculate the relative error of the mean estimate values and the true matrices
\begin{equation} \label{eq:relative.error}
	Err_K = \frac{\mathrm{norm}(K_+ -\bar{K}_{\mathrm{est}} )}{\mathrm{norm}(K_+)}  \qquad Err_R = \frac{\mathrm{norm}(R -\bar{R}_{\mathrm{est}} )}{\mathrm{norm}(R)}.
\end{equation}
Notice that $L$ might be different from the number of the sampled noisy $\Gamma$, which is 100 in our case, because of the possible failure in the optimization algorithm used in Algorithm~\ref{alg:reconstruction}. In addition, we test the approximation property of the trajectories generated by $(K_{\mathrm{est}},R_{\mathrm{est}})$ by comparing them with $\Gamma$ on plots. 

\begin{figure}
		\includegraphics[width=0.7\textwidth]{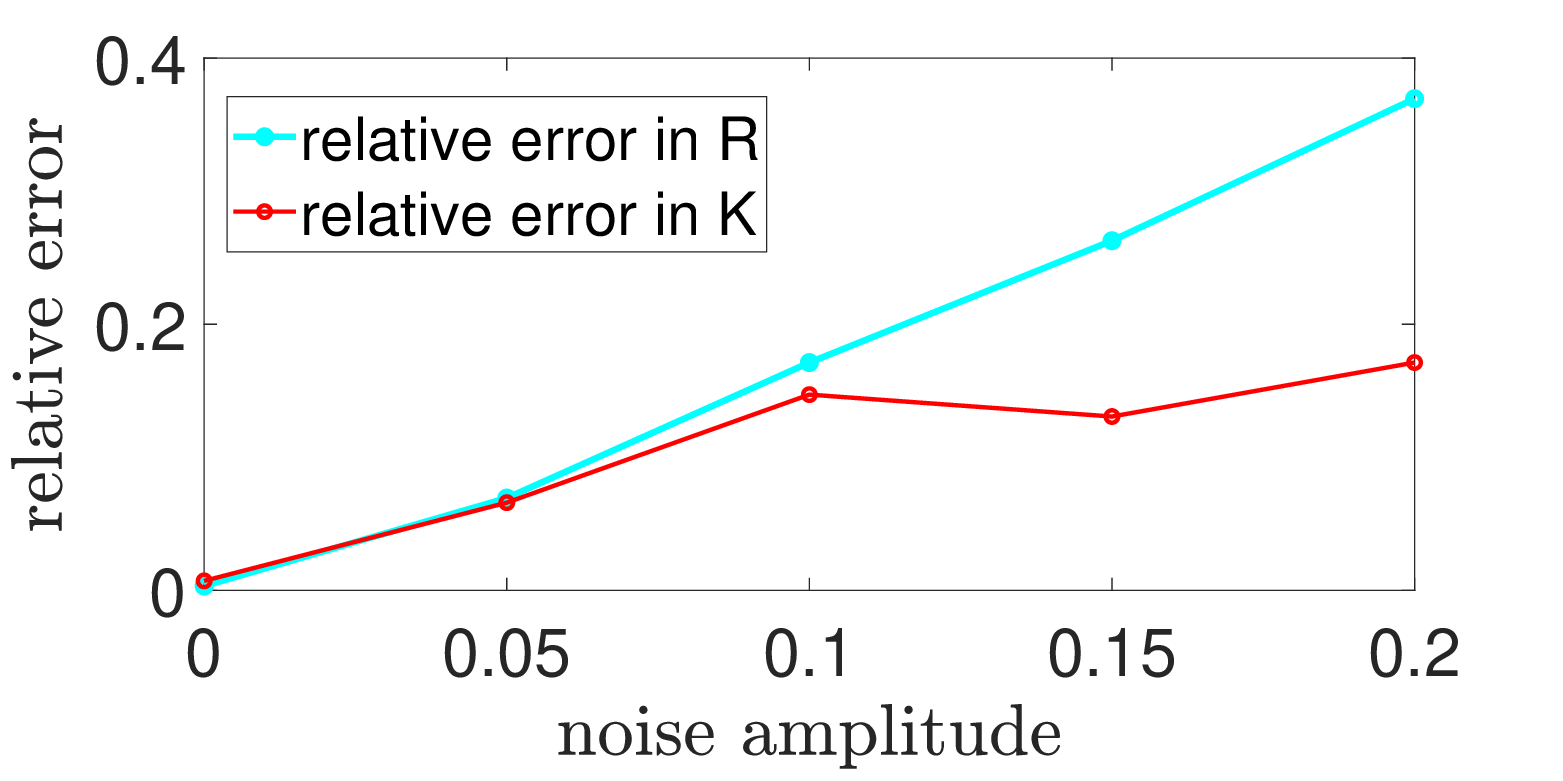}
	\caption{Plots of relative error for $(K_{\mathrm{est}},R_{\mathrm{est}})$ computed using~\eqref{eq:relative.error} for noise amplitudes $\alpha = 0, \ 0.05,\ 0.1,\ 0.15,\ 0.2$. }
\label{fig:relative.error}
\end{figure}

\begin{figure}
	\includegraphics[width=0.7\textwidth]{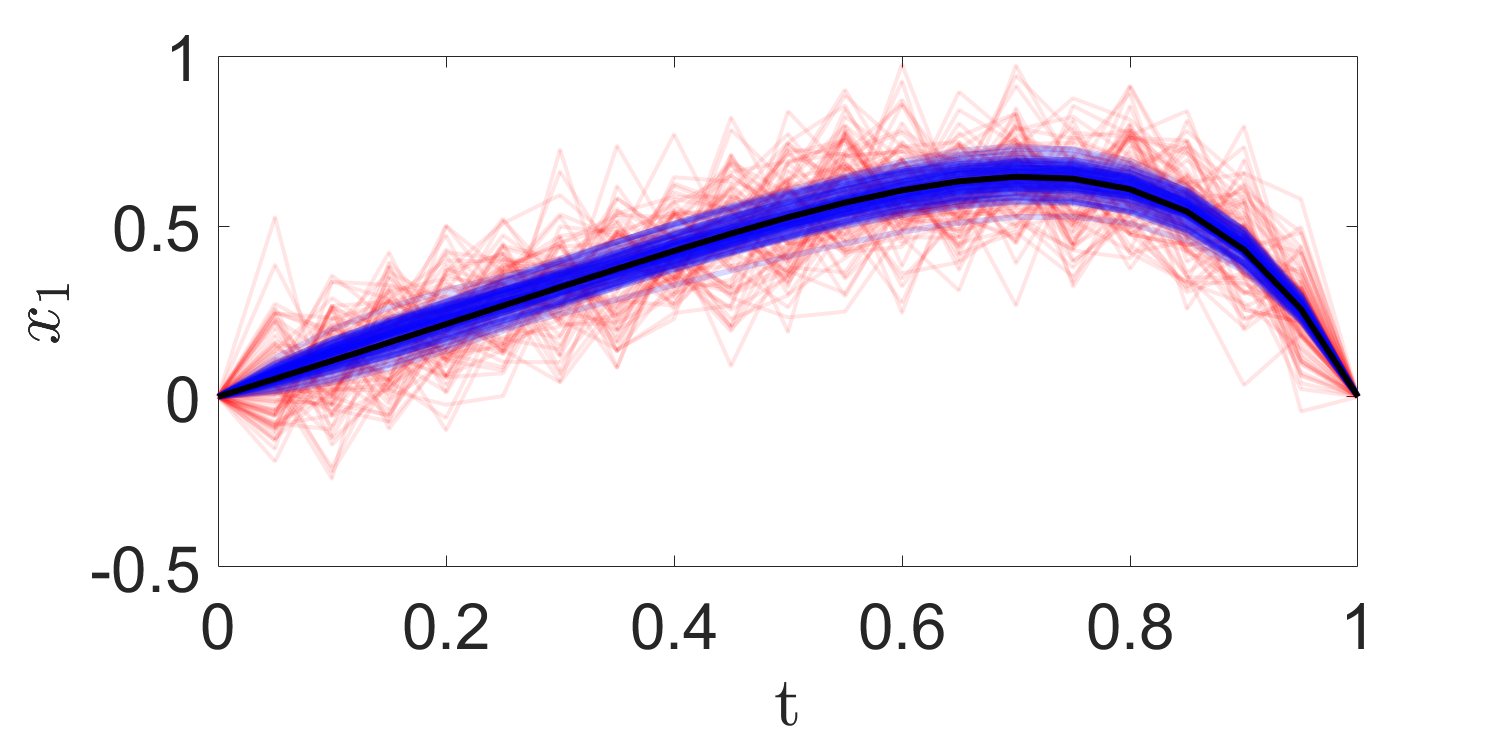}
\\[\smallskipamount]
	\includegraphics[width=0.7\textwidth]{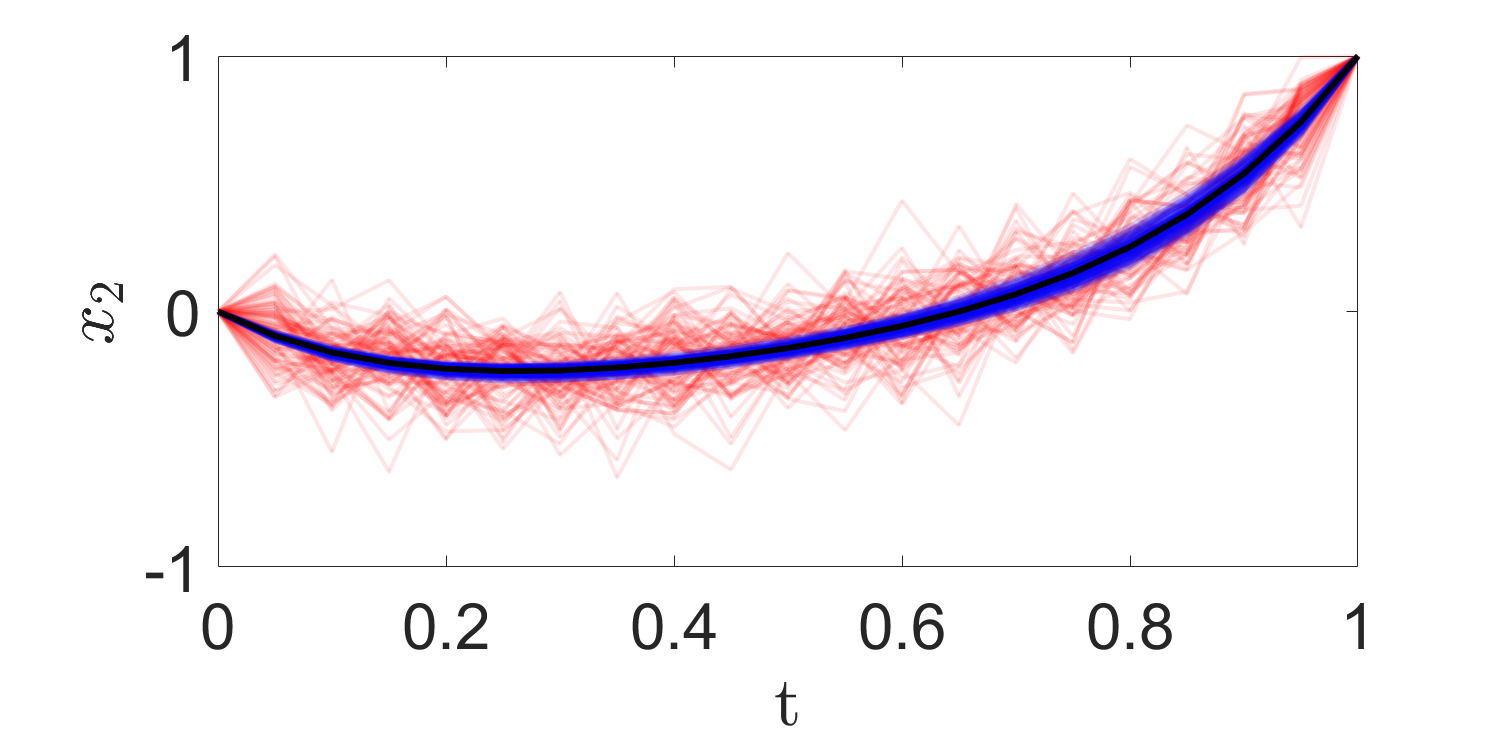}
\\[\smallskipamount]
	\includegraphics[width=0.7\textwidth]{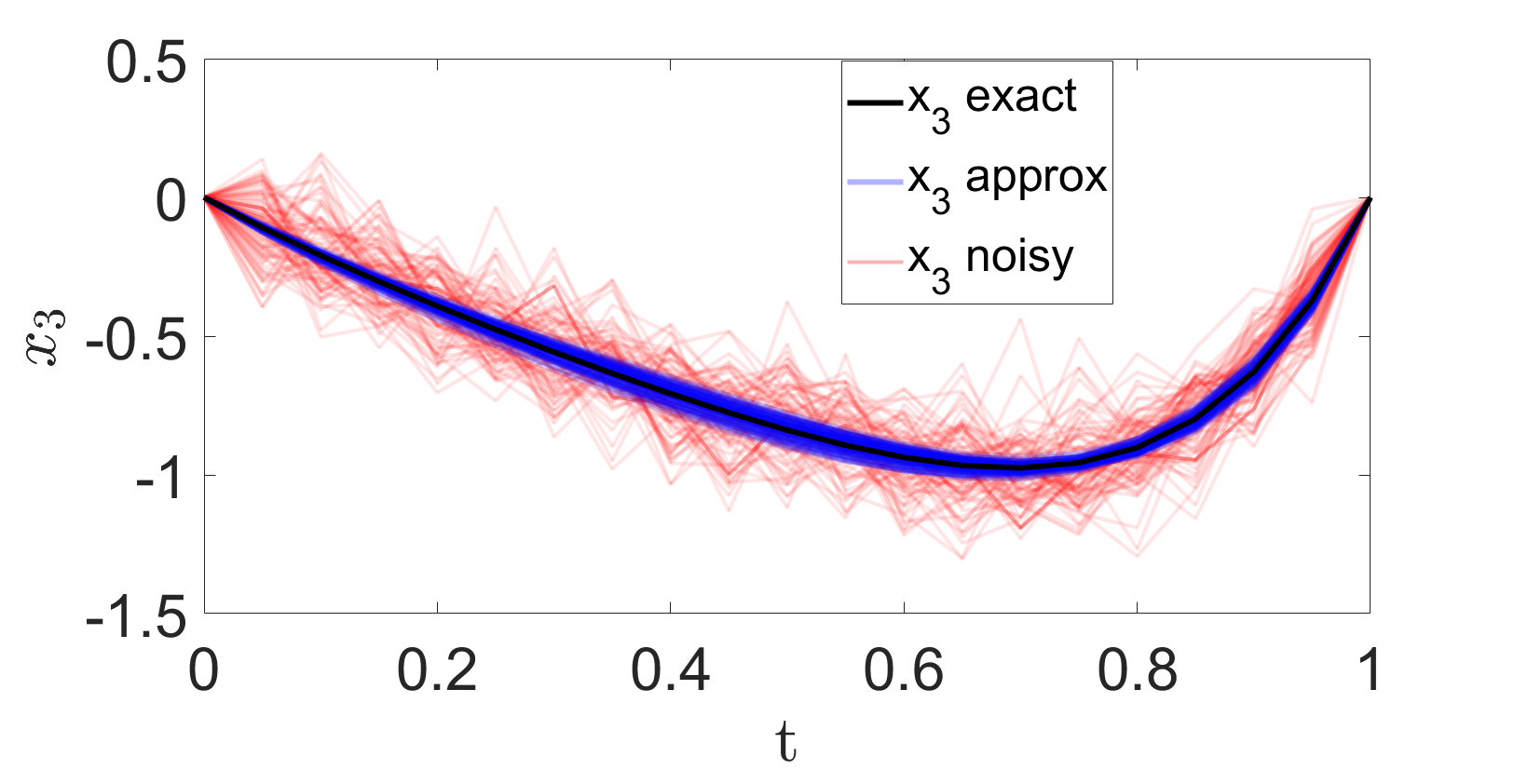}
\caption{Plots of the reconstructed trajectories in case of noise amplitude $\alpha = 0.2$. 
	The plots show the evolution of the coordinates $(x_1(t), x_2(t),x_3(t))$ separately. The red trajectories denote the noisy data obtained in $L=90$ samples of noise, the black curves are the exact solutions without noise, and the blue lines depict the trajectories obtained with the reconstructed cost from each sample of the noisy data.}
\label{fig:trajectories}
\end{figure}

In the example, we consider an LQ problem  with the control system given by
$$ A_1 =
\left(
\begin{array}{ccc}
	1 & 0 & 1 \\
	-2 & -3 & -1 \\
	0 & 0 & 2
\end{array}
\right), \quad B_1 = \left(
\begin{array}{cc}
	1 & 0 \\
	0 & 1 \\
	0 & 1
\end{array}
\right).$$
The cost is given by matrices $Q, \ S, \ R$ defined as follows
\begin{equation*}
	Q =
	\left(
	\begin{array}{ccc}
		20  &  6  &  34 \\
		6  &   2  &  11 \\
		34  &  11 &   61
	\end{array}
	\right), \quad S = \left(
	\begin{array}{cc}
		20 & 12 \\
		6   &  4 \\
		34  & 22
	\end{array}
	\right), \quad R = \left(
	\begin{array}{cc}
		5 & 3 \\
		3 & 2
	\end{array}
	\right).
\end{equation*}
%
%
%
%
Figures~\ref{fig:relative.error},~\ref{fig:trajectories} present the obtained results for different amplitudes of noise. Figure~\ref{fig:relative.error} shows the relative error of the estimated matrices $K, R$ calculated using \eqref{eq:relative.error}. Figure~\ref{fig:trajectories} presents the trajectories obtained using $K, R$ reconstructed from $L_{0.2}$ samples of noisy trajectories in case of $\alpha = 0.2$. Notice that the number of the successfully finished optimization algorithms and, thus, the number of the used reconstruction results out of 100 samples is $L_{0.0} = 100$, $L_{0.05} = 99$, $L_{0.1} = 98$, $L_{0.15} = 98$, $L_{0.2} = 90$.

We observe in Figure~\ref{fig:relative.error} that the relative error does not exceed 0.4 for all considered noise amplitudes. In Figure~\ref{fig:trajectories} we observe that the trajectories obtained from estimated pairs $(K_{\mathrm{est}},R_{\mathrm{est}})$ stay close to the true solution. We can conclude from these numerical results that  Algorithm~\ref{alg:reconstruction} provides stable numerical solutions for the inverse problem applied to noisy data and therefore, can be efficiently applied to measurements coming from applications.


\appendix

\section{Proof of Theorem~\ref{th:decomp_sol}}
\label{se:appendix} 

The proof of the theorem relies on several preliminary lemmas. The first one is a direct consequence of a classical result on Riccati equations.

\begin{Lemma}
\label{le:P+}
Under Assumptions~\ref{Assump:ABQSR},  \ref{cond:A1}, and \ref{cond:unif_contr}, there exist a  bounded and nonnegative solution $P_+ : \R \to \Sym(n)$ of the Riccati differential equation~\eqref{eq:Riccati_edo} such that:
\begin{itemize}
  \item $A_{+} = A-B R^{-1} \left(S^\top + B^\top P_{+}\right)$
is exponentially stable;
  \item for every $t_0 \in \R$ and $x_0 \in \R^n$, the feedback control defined by $u_+(t)=-R(t)^{-1} \left(S(t)^\top + B(t)^\top P_{+}(t)\right) x(t)$ minimizes
\begin{equation}
\lim_{t_1 \to +\infty} \left( C_{t_0}^{t_1}(x_u) + \|x_u(t_1)\|^2 \right), \qquad \hbox{where } x_u(t_0)=x_0,
\end{equation}
and the value of the minimum is
\begin{equation}
x_0^\top P_+(t_0) x_0.
\end{equation}
\end{itemize}
\end{Lemma}

\begin{proof}
Assume first $S=0$ and $A$ is exponentially stable. In that case, the lemma results from~\cite[Th.\ 3.5 and 3.6]{Kwakernaak1972} and~\cite[Th.\ 7.5.3]{AbouKandil2003} (the latter for the boundedness of $P_+$), which deal with costs of the form $c(t,x,u)=x^\top Q(t) x + u^\top R(t) u$. When $S \neq 0$, we write our quadratic cost as
\begin{equation}
c(t,x,u)=  x^{\top}Q_S(t) x + \bar{u}^{\top} R(t) \bar{u},
\end{equation}
where $Q_S = Q - S R^{-1} S^\top$, $\bar{u}=u +R^{-1} S^\top x$, and we write the control system as $\dot{x} = A_Sx + B \bar{u}$, where $A_S=A-BR^{-1} S^\top$. Applying the previous results, we obtain that the conclusions of the lemma hold if $A_S$ is exponentially stable.

Now, condition~\ref{cond:unif_contr} implies that $(A_S,B)$ is uniformly controllable, and then that there exists a bounded mapping $K : \R \to M_{m,n}(\R)$ such that $A_S-BK$ is exponentially stable. As above, lets us make a feedback transformation on the control. By setting $\tilde{u}=u+Kx$, we can write our quadratic cost as
\begin{equation}
\label{eq:cost_somme_carres}
c(t,x,u)=  x^{\top}\tilde{Q}(t) x + 2 x^{\top}\tilde{S}(t)\tilde{u} + \tilde{u}^{\top} R(t) \tilde{u},
\end{equation}
where $\tilde{Q} = Q - S K - K^\top S + K^\top R K$, $\tilde{S}=S-K^\top R$, and we write the control system as $\dot{x} = \tilde{A}x + B \tilde{u}$, where $\tilde{A}=A-BK$. It is easy to check that $(\tilde{A},B, \tilde{Q},\tilde{S},R)$ satisfy Assumptions  \ref{cond:A1}, \ref{cond:unif_contr}, and~\ref{Assump:ABQSR}, and since $\tilde{A}- BR^{-1}\tilde{S}^\top=A_S-BK$ is exponentially stable, the conclusions of the lemma hold for these cost and control system. We thus obtain matrices $\tilde{P}_+$, $\tilde{A}_+$, and a stabilizing control $\tilde{u}_+$.

A standard computation shows that the Riccati differential equations defined by $(\tilde{A},B, \tilde{Q},\tilde{S},R)$ and $(A,B,Q,S,R)$ are identical, so $\tilde{P}_+$ is a solution of~\eqref{eq:Riccati_edo}, and
\begin{equation}
\tilde{A}_{+} = \tilde{A}-B R^{-1} \left(\tilde{S}^\top + B^\top \tilde{P}_{+}\right)=A-B R^{-1} \left(S^\top + B^\top \tilde{P}_{+}\right).
\end{equation}
Lemma~\ref{le:P+} then follows by setting $P_+=\tilde{P}_+$.
\end{proof}


\begin{Lemma}
\label{le:P-}
Under Assumptions~\ref{Assump:ABQSR},  \ref{cond:A1}, and \ref{cond:unif_contr},  there exist a  bounded and  nonpositive solution $P_- : \R \to \Sym(n)$ of the Riccati differential equation~\eqref{eq:Riccati_edo} such that:
\begin{itemize}
  \item $A_{-} = A-B R^{-1} \left(S^\top + B^\top P_{-}\right)$
is exponentially anti-stable;
  \item for every $t_1 \in \R$ and $x_1 \in \R^n$, the feedback control defined by $u_-(t)=-R(t)^{-1} \left(S(t)^\top + B(t)^\top P_{-}(t)\right) x(t)$ minimizes
\begin{equation}
\lim_{t_0 \to -\infty} \left( C_{t_0}^{t_1}(x_u) + \|x_u(t_0)\|^2 \right), \qquad \hbox{where } x_u(t_1)=x_1,
\end{equation}
and the value of the minimum is
\begin{equation}
- x_1^\top P_-(t_1) x_1.
\end{equation}
\end{itemize}
\end{Lemma}

\begin{proof}
Fix $t_1\in \R$. For any trajectory $x$ of $\dot x = Ax + Bu$ associated with a control law $u$, the function $\tilde{x}(t)=x(2t_1-t)$ is  a trajectory associated with the control $\tilde{u}(t) =u(2t_1-t)$ of the control system
\begin{equation}
\dot{\tilde{x}}=\tilde{A} \tilde{x} + \tilde{B} \tilde{u},
\end{equation}
where the mappings $\tilde{A}, \tilde{B}$ are defined by $\tilde{A}=-A(2t_1-t)$, $\tilde{B}=-B(2t_1-t)$. Moreover, for $t_0<t_1$ we have
\begin{equation}
C_{t_0}^{t_1}(x_u) = \tilde{C}_{t_1}^{2t_1-t_0}(\tilde{u}),
\end{equation}
where the integral cost $\tilde{C}_{t_1}^{2t_1-t_0}$ is defined by means of the mappings $\tilde{Q}(t) =Q(2t_1-t)$, $\tilde{S} (t)=S(2t_1-t)$, $\tilde{R} (t) =R(2t_1-t)$. The mappings $(\tilde{A},\tilde{B},\tilde{Q},\tilde{S},\tilde{R})$ obviously satisfy Assumptions~\ref{Assump:ABQSR},  \ref{cond:A1}, and \ref{cond:unif_contr}, therefore by Lemma~\ref{le:P+} there exists a  bounded and nonnegative solution $\tilde{P}_+ : \R \to \Sym(n)$ of the Riccati differential equation associated with $(\tilde{A},\tilde{B},\tilde{Q},\tilde{S},\tilde{R})$. The conclusion follows by noticing that
\begin{equation}
P_-(t)=-\tilde{P}_+ (2t_1-t)
\end{equation}
is solution of the Riccati differential equation~\eqref{eq:Riccati_edo}, and that
\begin{equation}
A_{-}(t) = \tilde{A}_+(2t_1-t)
\end{equation}
is exponentially anti-stable and satisfies $A_-= A-B R^{-1} \left(S^\top + B^\top P_{-}\right)$.
\end{proof}

We are now in a position to prove the theorem.

\begin{proof}[Proof of Theorem~\ref{th:decomp_sol}]
Lemmas~\ref{le:P+} and \ref{le:P-} ensure the existence of the mappings $P_+,P_-$, the uniqueness resulting from \cite[Prop.\ 7.2.9]{AbouKandil2003}. Hence Points \emph{(i)-(ii)} are proved.
\smallskip

Let us prove \emph{(iii)}. We have to show that the nonnegative matrix $P_+(t)-P_-(t)$ is injective for every $t$.
Let $t_0 \in \R$ and $x_0 \in \R^n$ be such that $x_0 \in \ker (P_+(t_0)-P_-(t_0))$. Necessarily, $x_0$ belongs to the kernel of both $P_{\pm}(t_0)$, so that $x_0^\top P_{\pm}(t_0) x_0=0$. On the other hand, from the last point of Lemma~\ref{le:P+}, there holds
\begin{equation}
0=x_0^\top P_+(t_0) x_0 = \int_{t_0}^{+ \infty} c(t,x_+(t),u_+(t)) dt,
\end{equation}
where $x_+(\cdot)$ denotes the trajectory associated with the control $u_+$ such that $x_+(t_0)=x_0$. In particular $x_+$ is a solution of  $\dot{x}_+=A_+ x_+$ and $x_+(t) \to 0$ as $t \to + \infty$.  As in the proof of Lemma~\ref{le:P+}, we set $\bar{u}=u_+ +R^{-1} S^\top x_+$ and we obtain
\begin{equation}
0= \int_{t_0}^{+ \infty} c(t,x_+(t),u_+(t)) dt = \int_{t_0}^{+ \infty} \left(x_+(t)^{\top}Q_S(t) x_+(t) + \bar{u}(t)^{\top} R(t) \bar{u}(t) \right) dt.
\end{equation}
By \ref{cond:Rdefpos} and~\ref{cond:barQpos}, we get $\bar{u}(t)=0$  for every $t \geq t_0$. Hence $x_+$ satisfies $\dot{x}_+= Ax_+ - B R^{-1} S^\top x_+$ on $[t_0,+ \infty)$. In other terms, for every $t \geq t_0$, we have $x_+(t)=\Phi_{A_S}(t,t_0)x_0$, where $A_S= A - B R^{-1} S^\top$, and thus $\Phi_{A_S}(t,t_0)x_0  \to 0$ as $t \to + \infty$.

The same reasoning on $P_-$ shows that $\Phi_{A_S}(t,t_0)x_0  \to 0$ as $t \to - \infty$. However, arguing as in the proof of Lemma~\ref{le:P+}, we can assume up to a feedback transformation that the mapping $A_S$ is exponentially stable. And by definition of exponential stability, $\Phi_{A_S}(t,t_0)x_0  \to 0$ as $t \to \pm \infty$ implies $x_0=0$. Hence $\ker (P_+(t)-P_-(t))=\{0\}$ for every $t \in \R$ and \emph{(iii)} is proved.
\smallskip

It remains to prove \emph{(iv)}. First it is classical that Assumptions~\ref{Assump:ABQSR},  \ref{cond:A1}, and \ref{cond:unif_contr} imply existence and uniqueness of the optimal trajectory (remind that under theses hypothesis the cost written in the form~\eqref{eq:cost_somme_carres} is nonnegative and strictly convex w.r.t.\ the control). For given data $t_0,t_1,x_0,x_1$, this latter optimal trajectory is the first component $x^*$ of the unique solution $(x^*,p^*)$ of the Hamiltonian differential equation \eqref{eq:ham_edo} such that $x^*(t_0) = x_0$ and $x^*(t_1) = x_1$.

Let us study the set $\hh$ of all solutions of the Hamiltonian differential equation \eqref{eq:ham_edo}. The equation being linear, $\hh$ is a linear subspace of dimension $2n$ of the set of smooth mappings from $\R$ to $\R^{2n}$. Choose a symmetric solution $P$ of the Riccati differential equation~\eqref{eq:Riccati_edo} and a solution $x$ of the differential equation
\begin{equation}
\dot{x}= \left(A-B R^{-1}(S^\top + B^\top P)\right)x.
\end{equation}
It is just a matter of a calculation to check that the function $(x,-Px)$ belongs to $\hh$. We can choose for instance $P=P_+$ or $P_-$, and we obtain that the functions of the form
\begin{equation}
t \mapsto
\begin{pmatrix}
x(t) \\ p(t)
\end{pmatrix} =
\begin{pmatrix}
  I \\ -P_{\pm}(t)
\end{pmatrix} \Phi_{A_{\pm}}(t,t_0) y, \qquad y \in \R^n, \ t_0 \in \R,
\end{equation}
belong to $\hh$. Hence we can define a linear application $F_{t_0}: \R^{2n} \to \hh$ by
\begin{equation}
F_{t_0}(y_+,y_-)(t)= \begin{pmatrix}
  I \\ -P_{+}(t)
\end{pmatrix} \Phi_{A_{+}}(t,t_0) y_+ + \begin{pmatrix}
  I \\ -P_{-}(t)
\end{pmatrix} \Phi_{A_{-}}(t,t_0) y_-, \quad \forall t \in \R.
\end{equation}
This application is injective. Indeed, if $(y_+,y_-)$ belongs to $\ker F_{t_0}$, then
\begin{equation}
0= F_{t_0}(y_+,y_-)(t_0)=
\begin{pmatrix}
  y_+ + y_- \\ -P_+ (t_0) y_+ - P_-(t_0) y_-,
\end{pmatrix}
\end{equation}
which implies $y_+=-y_-$ and $(P_+ (t_0) - P_-(t_0)) y_+=0$. By \emph{(iii)}, $P_+ (t_0) - P_-(t_0)$ is positive definite, which implies that $y_+=0$, $y_-=0$, and thus that $F_{t_0}$ is injective.

Since $\dim \hh=2n$, we conclude that $F_{t_0}$ is an isomorphism  between $\R^{2n}$ and $\hh$. As a consequence, there exists a unique pair $(y_+,y_-)$ such that $(x^*,p^*)=F_{t_0}(y_+,y_-)$, which proves point \emph{(iv)} of the theorem.
\end{proof}

\section{Linearisation of Hamiltonian flows}

Let us consider a (nonlinear) control system of the form
\begin{equation}
	\label{eq:dymsys_nl}
	\dot{x} = f(x,u)= f_0(x) + \sum_{i=1}^m u_i f_i(x), \qquad x\in \R^n, \quad u \in \R^m,
\end{equation}
where $f_0,\dots,f_m$ are vector fields on $\R^n$, and two quadratic cost
\begin{equation}
	\label{eq:qvcost_nl}
	c_i(x,u)=   u^{\top} R_i(x) u, \quad \hbox{where } R_i(x) \succ 0 \quad \forall x \in \R^n, \qquad i=1,2.
\end{equation}
For $i=1,2$, the associated normal Hamiltonian $h_i: T^*\R^n=\R^n \times (\R^n)^* \to \R$ are defined as $h_i(x,p)=\max \{ \widetilde{h}_i(x,p,u,- 1/2) \ : \ u \in \R^m \}$ where $\widetilde{h}_i(x,p,u,p^0) = \langle p, f(x,u) \rangle + p^0 c_i(x,u)$. A normal extremal of $c_i$ is a trajectory $(x,p)(\cdot)$ of the Hamiltonian vector field $\vec{h}_i$, i.e.,
\begin{equation}
	(x,p)(t) = e^{t \vec{h}_i} ((x,p)(0)).
\end{equation}
A normal geodesic of $c_i$ is the projection $x(\cdot)$ on $\R^n$ of a normal extremal, i.e., there exists $p_0 \in (\R^n)^*$ such that $x(t)=\pi \left( e^{t \vec{h}_i} (x(0),p_0)\right)$, where $\pi$ is the canonical projection $\pi (x,p)=x$.

Suppose that the optimal control problems associated with \eqref{eq:dymsys_nl}, $c_1$ and \eqref{eq:dymsys_nl}, $c_2$ have the same optimal synthesis, namely that any trajectory of \eqref{eq:dymsys_nl} minimizing the integral cost $\int c_1(x,u)$ between its extremities also minimizes $\int c_2(x,u)$, and vice versa. In other terms, $c_1$ and $c_2$ are two solutions of the same inverse optimal control problem.

Consider a minimizer $(\bar{x}(t),\bar{u}(t))$, $t \in \R$, of $c_1$ such that the matrix-valued mapping $(A,B)$ defined by
\begin{equation}\label{eq:linearised}
	A(t) = \frac{\partial f}{\partial x}\left(\bar{x}(t),\bar{u}(t) \right), \quad B(t) = \frac{\partial f}{\partial u}\left(\bar{x}(t),\bar{u}(t) \right), \qquad t \in \R,
\end{equation}
satisfy the Kalman sufficient condition~\ref{cond:kalman}. By Pontryagin Maximum Principle and \cite[Prop.\ 3.12]{AgrachevBarilariRizzi}, such a minimizer is necessarily a normal geodesic of $c_1$, therefore there exists $p_0 \in (\R^n)^*$ such that $\bar{x}(t)=\pi \left( e^{t \vec{h}_1} (\bar{x}_0,p^1_0)\right)$ with $\bar{x}_0=\bar{x}(0)$. The covector $p^1_0$ is said to be \emph{ample with respect to $c_1$}. Then, from \cite[Prop.\ 3.6]{Jean2019a}, there exists open subset $V_1 \ni (\bar{x}_0,p^1_0)$ and $V_2$ of  $T^*\R^n$ and a diffeomorphism $\phi: V_1 \rightarrow V_2$ such that $\pi \circ \phi = \pi$ and $\phi$ sends the integral curves of $\vec h_1$ to the integral curves of $\vec h_2$, i.e.\ $\phi \bigl(e^{t\vec h_1}(x,p)\bigr) =e^{t\vec h_2} \bigl(\phi(x,p)\bigr)$ for all $(x,p) \in V_1$ and $t\in \mathbb R$.

By differentiation of the latter relation at $(\bar{x}_0,p^1_0)$, we obtain
\begin{equation}\label{eq:diff_orbital_diffeo}
	D\phi \bigl(e^{t\vec h_1}(\bar{x}_0,p^1_0)\bigr) \circ De^{t\vec h_1}(\bar{x}_0,p^1_0)  = D e^{t\vec h_2} \bigl(\phi(\bar{x}_0,p^1_0)\bigr) \circ D\phi(\bar{x}_0,p^1_0).
\end{equation}

A standard computation shows that the differential of the flow of the Hamiltonian vector fields $\vec{h}_i$ is the flow of a linear Hamiltonian vector field arising from a LQ problem. More precisely, $De^{t\vec h_1}(\bar{x}_0,p^1_0) = \Phi_{\vec{H}_1}(t,0)$, where $\Phi_{\vec{H}_1}$ is the transition matrix of $\vec{H}_1(\cdot)$ and, for every $t \in \R$, $\vec{H}_1(t)$ is the Hamiltonian matrix of the form~\eqref{eq:Ham_matrix} defined by the matrices $A(t),B(t)$ given in~\eqref{eq:linearised} and
\begin{multline}
	Q_1(t)= \frac{\partial^2 \widetilde{h}_1}{\partial x^2} \left(\bar{x}(t),\bar{p}^1(t),\bar{u}(t),- \frac{1}{2} \right), \quad
	S_1(t)= \frac{\partial^2 \widetilde{h}_1}{\partial x \partial u} \left(\bar{x}(t),\bar{p}^1(t),\bar{u}(t),- \frac{1}{2} \right),\\
	R_1(t)= \frac{\partial^2 \widetilde{h}_1}{\partial u^2} \left(\bar{x}(t),\bar{p}^1(t),\bar{u}(t),- \frac{1}{2} \right), \qquad \qquad
\end{multline}
the covector $\bar{p}^1(t)$ being defined by $(\bar{x}(t),\bar{p}^1(t)) = e^{t\vec h_1}(\bar{x}_0,p^1_0)$. In the same way, $D e^{t\vec h_2} \bigl(\phi(\bar{x}_0,p^1_0)\bigr)= \Phi_{\vec{H}_2}(t,0)$, where, for every $t \in \R$, $\vec{H}_2(t)$ is the Hamiltonian matrix of the form~\eqref{eq:Ham_matrix} defined by the matrices $A(t),B(t)$ and
\begin{multline}
	Q_2(t)= \frac{\partial^2 \widetilde{h}_2}{\partial x^2} \left(\bar{x}(t),\bar{p}^2(t),\bar{u}(t),- \frac{1}{2} \right), \quad
	S_2(t)= \frac{\partial^2 \widetilde{h}_2}{\partial x \partial u} \left(\bar{x}(t),\bar{p}^2(t),\bar{u}(t),- \frac{1}{2} \right),\\
	R_2(t)= \frac{\partial^2 \widetilde{h}_2}{\partial u^2} \left(\bar{x}(t),\bar{p}^2(t),\bar{u}(t),- \frac{1}{2} \right), \qquad \qquad
\end{multline}
the covector $\bar{p}^2(t)$ being defined by $(\bar{x}(t),\bar{p}^2(t)) = e^{t\vec h_2}\bigl(\phi(\bar{x}_0,p^1_0)\bigr)$ (notice that $(\bar{x}(t),\bar{p}^2(t)) =\phi (\bar{x}(t),\bar{p}^1(t))$.

Set $\psi_t= D\phi \bigl(e^{t\vec h_1}(\bar{x}_0,p^1_0)\bigr)$. For every $t \in \R$, $\psi_t : T^*\R^n \to T^*\R^n$ is a linear isomorphism and satisfies $\pi \circ \psi_t = \pi$. With this notation, \eqref{eq:diff_orbital_diffeo} writes as
\begin{equation}
	\psi_t \circ \Phi_{\vec{H}_1}(t,0) = \Phi_{\vec{H}_2}(t,0) \circ \psi_0.
\end{equation}
As a consequence, for every $(\delta x_0,\delta p_0) \in T^*\R^n$, we have
\begin{equation}
	\pi \circ \Phi_{\vec{H}_1}(t,0) (\delta x_0,\delta p_0)  = \pi \circ \Phi_{\vec{H}_2}(t,0) (  \psi_0 (\delta x_0,\delta p_0)).
\end{equation}
This equality exactly means that all the minimizing solutions of the LQ problem defined by $(A,B,Q_1,S_1,R_1)$, which can be written as $\pi \circ \Phi_{\vec{H}_1}(t,0) (\delta x_0,\delta p_0)$ by \eqref{eq:ham_edo}, are also minimizing solutions of the LQ problem defined by $(A,B,Q_2,S_2,R_2)$. The map $\psi_0$ being an isomorphism, we obtain that  $(A,B,Q_1,S_1,R_1)$ and $(A,B,Q_2,S_2,R_2)$ define the same LQ optimal synthesis.

To summarize, fix a control-affine system~\eqref{eq:dymsys_nl} and assume that $c_1$, $c_2$ are solutions of the same inverse optimal control problem. Then, the linearisation of the optimal control problem around any ample geodesic gives rise to two cost solutions of the same LQ inverse optimal control problem. This is one of the main motivation for the study of time-dependent LQ inverse optimal control problems.

\bibliographystyle{plain}
\bibliography{Biblio/biblio_LQ}

\end{document}